\pdfoutput=1
\documentclass[11pt]{amsart}
\usepackage[T1]{fontenc}

\usepackage[margin=1in]{geometry}

\usepackage{amsmath}
\usepackage{amsthm}
\usepackage{mathtools}
\usepackage{amssymb}

\usepackage{newtxtext}
\usepackage{newtxmath}
\usepackage[cal=euler]{mathalpha}

\usepackage[round]{natbib}
\usepackage{tikz}
\usetikzlibrary{arrows.meta,calc,positioning}
\usepackage{enumitem}
\usepackage[hidelinks]{hyperref}
\usepackage[capitalise,noabbrev]{cleveref}
\newcommand{\doi}[1]{\url{https://doi.org/#1}}

\crefname{mainthm}{Theorem}{Theorems}
\Crefname{mainthm}{Theorem}{Theorems}
\crefname{prop}{Proposition}{Propositions}
\Crefname{prop}{Proposition}{Propositions}
\crefname{examples}{Example}{Examples}
\Crefname{examples}{Example}{Examples}
\crefname{notation}{Notation}{Notations}
\Crefname{notation}{Notation}{Notations}
\usepackage{microtype}
\usepackage{stmaryrd}
\usepackage{aliascnt}

\newcommand{\sharedtheorem}[2]{%
  \newaliascnt{#1}{theorem}%
  \newtheorem{#1}[#1]{#2}%
  \aliascntresetthe{#1}}
\newtheorem{theorem}{Theorem}[section]
\newtheorem{mainthm}{Theorem}

\sharedtheorem{corollary}{Corollary}
\sharedtheorem{prop}{Proposition}
\sharedtheorem{lemma}{Lemma}
\theoremstyle{definition}
\sharedtheorem{definition}{Definition}
\sharedtheorem{example}{Example}
\sharedtheorem{examples}{Examples}
\sharedtheorem{notation}{Notation}
\theoremstyle{remark}
\sharedtheorem{remark}{Remark}
\newcommand{\cat}[1]{\mathbf{#1}}
\newcommand{\bool}{\cat{Bool}}
\newcommand{\law}{\cat{Cost}}
\newcommand{\elaw}{\overline{\cat{Cost}}}
\newcommand{\C}{{\underline{\cat{C}}}}
\newcommand{\terminal}{\underline{\cat{1}}}
\renewcommand{\S}{\cat{S}}
\newcommand{\D}{{\underline{\cat{D}}}}
\newcommand{\E}{{\underline{\cat{E}}}}
\DeclareMathOperator{\obj}{ob}
\newcommand{\Cat}{\cat{Cat}}

\newcommand{\Adj}{\cat{Adj}}
\newcommand{\quant}{\cat{Q}}
\newcommand{\Quant}{{\underline{\cat{Q}}}}
\newcommand{\QCat}{\quant\Cat}
\newcommand{\QAdj}{\quant\Adj}
\newcommand{\op}{\mathrm{op}}
\newcommand{\meet}{\wedge}
\newcommand{\bigmeet}{\bigwedge}
\newcommand{\bigjoin}{\bigvee}
\newcommand{\dlhom}[2]{\{ #1, #2 \}}
\let\origcdot\cdot
\newcommand{\scdot}{\origcdot}
\renewcommand{\cdot}{\,\smash{\raisebox{0.08em}{$\scriptscriptstyle\bullet$}}\,}
\newcommand{\tensor}{\otimes}
\newcommand{\cotensor}{\pitchfork}
\DeclareMathOperator{\id}{id}
\DeclareMathOperator{\pre}{pre}
\DeclarePairedDelimiter{\norm}{\lVert}{\rVert}
\newcommand{\qequals}{\simeq}
\newcommand{\qleq}{\precsim}
\newcommand{\face}{\trianglelefteqslant}
\newcommand{\shf}{F}            
\newcommand{\ushf}{\underline{\shf}}      
\newcommand{\oshf}{{\overline{\shf}}}       
\newcommand{\dshf}{\widetilde{\shf}}       
\newcommand{\cochains}{\cat{C}^0(G;\shf)}
\newcommand{\wsections}{\cat{\Gamma}^W(G;\shf)}
\newcommand{\transport}[2]{\mathcal{P}^{\shf}_{#1 \to #2}}
\renewcommand{\preceq}{\preccurlyeq}
\renewcommand{\succeq}{\succcurlyeq}

\begin{document}

\title{Categorical Diffusion of Weighted Lattices}
\author[Robert Ghrist]{Robert Ghrist\textsuperscript{1}}
\author[Miguel Lopez]{Miguel Lopez\textsuperscript{2}}
\author[Paige Randall North]{Paige Randall North\textsuperscript{3}}
\author[Hans Riess]{Hans Riess\textsuperscript{4,\ensuremath{\ast}}}
\thanks{\textsuperscript{1}Department of Mathematics and Department of Electrical \& Systems Engineering, University of Pennsylvania}
\thanks{\textsuperscript{2}PlantingSpace}
\thanks{\textsuperscript{3}Department of Mathematics and Department of Information \& Computing Sciences, Utrecht University}
\thanks{\textsuperscript{4}School of Electrical \& Computer Engineering, Georgia Institute of Technology}
\thanks{\textsuperscript{\ensuremath{\ast}}Corresponding author: \href{mailto:riess@gatech.edu}{\texttt{riess@gatech.edu}}}
\date{}
\hypersetup{pdftitle={Categorical Diffusion of Weighted Lattices}, pdfauthor={Robert Ghrist, Miguel Lopez, Paige Randall North, Hans Riess}}

\begin{abstract}
Quantales are well known in logic and computer science, and categories enriched in them unify ordered, metric, relational, and tropical structures. Problems in dynamic programming, synchronization, and discrete event systems require reconciling local data, represented by objects in a quantale-enriched category, across the edges of a network, represented by an undirected graph. Concurrently, replacing inner products with hom-objects and linear adjoints with adjoint functors initiates a categorification of graph and sheaf Laplacians, with a view toward enrichment in larger monoidal categories and bicategories. To this end, we lay the foundations for an enriched harmonic analysis on graphs through cellular sheaves of quantale-enriched categories and a Laplace operator that we call the Lawvere Laplacian. Over an arbitrary unital quantale, the weighted global sections of such a sheaf are the vertex data that agree across every edge up to the tolerances of a weight function. They form the quantale-enriched category of natural transformations from a weight diagram into the sheaf, which is a weighted limit when the quantale is commutative. When the restriction maps have right adjoints, we prove that the weighted global sections are the suffix points of the Lawvere Laplacian, which are the equilibria of its harmonic flow. Under a continuity hypothesis, the harmonic flow converges to the greatest weighted global section below its initial cochain, and for central weight functions the weighted global sections form a complete coreflective subcategory of the cochains. Instances of the theory include the Bellman--Ford iteration over the Lawvere quantale, the reconciliation of local cost estimates, quantale and clock synchronization, and coupled max-plus discrete event systems.

\end{abstract}

\maketitle


\section{Introduction}\label{sec:intro}

On a graph endowed with a cellular sheaf of Hilbert spaces~\citep{hansen2019}, the degree-zero sheaf Laplacian aggregates local disagreement across edges: along each incidence it pulls data back from the edge stalk to the vertex stalk via the Hilbert-space adjoint of the restriction map, and the equilibria of its heat flow are the global sections of the sheaf.
This link between local diffusion and global consistency relies on an inner product to furnish adjoints and on the spectral theorem to guarantee convergence.

Many networked systems carry data that admit no inner product or additive inverse, such as truth values~\citep{fagiolini2008logical,fagiolini2016convergence}, path costs~\citep{mohri2002semiring,baras2010path}, and event timings~\citep{cohen1985linear,baccelli1992synchronization}.
Values of this kind form quantales, complete lattices endowed with an associative multiplication that distributes over joins: disjunction and conjunction for truth values, minimum and addition for costs, and maximum and addition for event timings.
Because aggregation proceeds by meets and joins, which are idempotent and noninvertible, linear Laplacians do not apply.
In this paper we replace Hilbert spaces with categories enriched in a unital quantale~\citep{Kelly2005,stubbe_introduction_2014}, linear adjoints with adjoint enriched functors, and sums with weighted meets.
The following problems motivate this paper:

\subsubsection*{Discrete event systems.}
In a production line or railway timetable, event epochs evolve by maxima of sums, making each subsystem linear over the max-plus semiring~\citep{cohen1985linear,baccelli1992synchronization,butkovic2010}, whose completion is isomorphic by negation to the extended Lawvere quantale (\Cref{ex:quantales}(e)). Classical max-plus models assemble coupled subsystems into a single matrix~\citep{heidergott2006maxplus,goverde2007railway}, obscuring the network decomposition. Modeling the subsystems as vertex stalks and their interactions along edges casts timetable synthesis as a local-to-global problem, whose harmonic flow computes the latest schedule below an initial bound that satisfies all synchronization tolerances (\Cref{ex:maxplus}).

\subsubsection*{Dynamic programming.}
Fixed-point iterations of min-plus Bellman operators solve the algebraic path problem~\citep{lehmann1977algebraic,baras2010path} and underlie dynamic programming and reinforcement learning~\citep{bertsekas2019reinforcement}. Over the Lawvere quantale, the Bellman--Ford shortest-path iteration is the harmonic flow of our Laplacian (\Cref{ex:shortest-paths}), and our convergence results recover its limit while generalizing to structured vertex stalks and costs in arbitrary quantales.

\medskip We lay the foundations for an enriched harmonic analysis on graphs through cellular sheaves of $\quant$-categories, governed by the \emph{Lawvere Laplacian}.
Its coefficients lie in a unital quantale $\quant$, and its data at each vertex form a $\quant$-category.
\citet{Law73} observed that $\quant$-categories unify preorders and metric spaces, which arise over the Boolean quantale and over the Lawvere quantale, respectively.
Complete $\quant$-categories, those admitting all weighted meets and joins, are the quantale-enriched analogues of complete lattices, and we call them \emph{weighted lattices}. Over $\quant = \bool$ they are the complete lattices, up to isomorphism of objects, and over $\quant = \law$ they are the Lawvere metric spaces that admit all weighted meets and joins, such as $[0,\infty]$ with truncated subtraction.
The data are organized by a cellular sheaf in the sense of \citet{shepardCellularDescriptionDerived1985} and \citet{curry2014}: a functor out of the incidence poset of the graph, which assigns data to the vertices and edges and maps to the incidences between them.
Such a functor determines a sheaf on the Alexandrov topology of the incidence poset (\Cref{rem:alexandrov}).
Generalizing the Tarski Laplacian of \citet{Ghrist_2022}, which acts on sheaves of complete lattices, the Lawvere Laplacian acts on sheaves of complete $\quant$-categories.
\Cref{fig:bridge} compares the resulting framework with classical harmonic analysis on graphs: the analogy with Hodge theory operates at degree zero, where the suffix points of the Lawvere Laplacian (the equilibria of the harmonic flow) coincide with the weighted global sections, just as the kernel of the sheaf Laplacian $L$ coincides with $H^0(G;\shf)$.

\begin{figure}[tp]
  \centering
  \begin{tikzpicture}[
      x=1cm, y=1cm,
      head/.style={align=center, font=\scriptsize\bfseries, text=black, inner sep=1pt},
      cpt/.style={align=center, font=\scriptsize\bfseries, text=black, inner sep=1pt},
      cls/.style={align=center, font=\scriptsize, text=black, inner sep=1pt},
      catc/.style={align=center, font=\scriptsize, text=black, inner sep=1pt},
      refc/.style={align=center, font=\scriptsize, text=black!70, inner sep=1pt},
      desc/.style={align=center, font=\scriptsize, text=black!75, inner sep=1pt},
      rline/.style={black!28, line width=0.4pt},
    ]

    \begin{scope}
      \clip[rounded corners=3pt] (0.00,0.00) rectangle (14.45,10.95);

      \fill[black!7]  (0.00,0.00)  rectangle (1.85,9.72);
      \fill[black!2]  (1.85,0.00)  rectangle (6.65,9.72);
      \fill[black!9]  (6.65,0.00)  rectangle (11.55,9.72);
      \fill[black!4]  (11.55,0.00) rectangle (14.45,9.72);

      \fill[black!15] (0.00,9.72)  rectangle (1.85,10.95);
      \fill[black!10] (1.85,9.72)  rectangle (6.65,10.95);
      \fill[black!18] (6.65,9.72)  rectangle (11.55,10.95);
      \fill[black!12] (11.55,9.72) rectangle (14.45,10.95);

      \foreach \xx in {1.85, 6.65, 11.55} {
        \draw[black!22, line width=0.35pt] (\xx,0.00) -- (\xx,10.95);
      }

      \draw[black!60, line width=0.65pt] (0.00,9.72) -- (16.00,9.72);
      \foreach \yy in {9.02, 8.00, 6.96, 5.98, 5.22, 4.24, 3.14, 2.10, 1.06} {
        \draw[rline] (0.00,\yy) -- (16.00,\yy);
      }
    \end{scope}

    \draw[rounded corners=3pt, black!55, line width=0.55pt]
      (0.00,0.00) rectangle (14.45,10.95);

    \node[head] at (0.925,10.33) {CONCEPT};
    \node[head] at (4.250,10.33)
      {SPECTRAL\\[0.5pt] SHEAF THEORY\\[1pt]
       {\normalfont\scriptsize\textcolor{black!65}{\citep{hansen2019}}}};
    \node[head] at (9.100,10.33)
      {QUANTALE-ENRICHED\\[0.5pt]
       DIFFUSION\\[1pt]
       {\normalfont\scriptsize\textcolor{black!65}{(this paper)}}};
    \node[head] at (13.000,10.33) {REFERENCE};

    \node[cpt]  at (0.925,9.37) {Coefficients};
    \node[cls]  at (4.250,9.37) {$\mathbb{R},\ \mathbb{C}$};
    \node[catc] at (9.100,9.37) {Unital quantale $\quant$};
    \node[refc] at (13.000,9.37) {\Cref{def:quantale}};

    \node[cpt]  at (0.925,8.51) {State space\\[-0.5pt]{\normalfont\scriptsize(cochains)}};
    \node[cls]  at (4.250,8.51)
      {Hilbert spaces\\[1pt]
       $C^{0}(G;\shf) = \bigoplus_{v \in V}\shf(v)$};
    \node[catc] at (9.100,8.51)
      {$\quant$-categories\\[1pt]
       $\cochains = \prod_{v \in V}\shf(v)$};
    \node[refc] at (13.000,8.51)
      {\Cref{def:q-category-main},\\[-0.5pt]\Cref{def:cochains}};

    \node[cpt]  at (0.925,7.48) {Similarity};
    \node[cls]  at (4.250,7.48)
      {Inner product\\[1pt]
       $\langle x,y\rangle = \sum_{v \in V}\langle x_v,y_v\rangle_{\shf(v)}$};
    \node[catc] at (9.100,7.48)
      {Hom-object\\[1pt]
       $\hom_{\cat{C}^{0}}(x,y) = \bigmeet_{v \in V}\hom_{\shf(v)}(x_v,y_v)$};
    \node[refc] at (13.000,7.48) {\Cref{def:cochains}};

    \node[cpt]  at (0.925,6.47) {Duality\\[-0.5pt]{\normalfont\scriptsize(adjoints)}};
    \node[cls]  at (4.250,6.47)
      {Linear operator and\\[0.5pt]adjoint $\shf_{v \face e}^{\ast}$};
    \node[catc] at (9.100,6.47)
      {Adjoint $\quant$-functors\\[1pt]
       $\ushf_{v \face e} \dashv_1 \oshf_{v \face e}$};
    \node[refc] at (13.000,6.47) {\Cref{def:q-adjoint-bisheaf}};

    \node[cpt]  at (0.925,5.60) {Transport};
    \node[cls]  at (4.250,5.60)
      {$\transport{w}{v} = \shf_{v \face e}^{\ast} \circ \shf_{w \face e}$};
    \node[catc] at (9.100,5.60)
      {$\transport{w}{v} = \oshf_{v \face e} \circ \ushf_{w \face e}$};
    \node[refc] at (13.000,5.60) {\Cref{def:transport}};

    \node[cpt]  at (0.925,4.73) {Aggregation};
    \node[cls]  at (4.250,4.73)
      {Linear sum\\[1pt]$\sum\nolimits_{w \sim_e v}$};
    \node[catc] at (9.100,4.73)
      {Weighted meet\\[1pt]$\bigmeet\nolimits_{w \sim_e v}^{W(v,-)}$};
    \node[refc] at (13.000,4.73) {\Cref{def:weighted}};

    \node[cpt]  at (0.925,3.69) {Laplacian};
    \node[cls]  at (4.250,3.69)
      {Sheaf Laplacian $(Lx)_v$\\[1.5pt]
       $\sum_{w \sim_e v}\bigl(\shf_{v \face e}^{\ast}\shf_{v \face e}x_v
         - \transport{w}{v}x_w\bigr)$};
    \node[catc] at (9.100,3.69)
      {Lawvere Laplacian $(L^{W}x)_v$\\[1.5pt]
       $\bigmeet\nolimits_{w \sim_e v}^{W(v,-)} \transport{w}{v}(x_w)$};
    \node[refc] at (13.000,3.69) {\Cref{def:lap}};

    \node[cpt]  at (0.925,2.62) {Dynamics};
    \node[cls]  at (4.250,2.62)
      {Heat semigroup $e^{-tL}$\\[1pt]
       $\partial_t x = -Lx$};
    \node[catc] at (9.100,2.62)
      {Harmonic flow $\Phi^{t}(x[0])$\\[1pt]
       $x[t+1]_v = x[t]_v \meet (L^{W}x[t])_v$};
    \node[refc] at (13.000,2.62) {\Cref{def:flow-functor}};

    \node[cpt]  at (0.925,1.58) {Equilibria};
    \node[cls]  at (4.250,1.58)
      {Harmonic cochains\\[1pt]
       $x(\infty) \in \ker L$};
    \node[catc] at (9.100,1.58)
      {Weighted global sections\\[1pt]
       $x[\infty] \in \wsections$};
    \node[refc] at (13.000,1.58)
      {\Cref{def:w-global-sections},\\[-0.5pt]\Cref{prop:kleene}};

    \node[cpt]  at (0.925,0.53) {Hodge\\[-0.5pt]theorem};
    \node[cls]  at (4.250,0.53)
      {Cellular Hodge theorem\\[1pt]
       $\ker L \cong H^{0}(G;\shf)$};
    \node[catc] at (9.100,0.53)
      {Hodge--Lawvere theorem\\[1pt]
       $\cat{Suffix}(L^{W}) = \wsections$};
    \node[refc] at (13.000,0.53) {\Cref{thm:hodge-thm}};

  \end{tikzpicture}
  \caption{Comparison between spectral sheaf theory~\citep{hansen2019} and the quantale-enriched diffusion theory of this paper. Cochains carry an inner product or hom-objects, transport composes a restriction with an adjoint, and neighborhood aggregation is a linear sum or a weighted meet. The heat flow converges to $\ker L \cong H^0(G;\shf)$ by the spectral theorem, and the harmonic flow converges, under the continuity hypothesis of \Cref{prop:kleene}, to the weighted global sections, which the Hodge--Lawvere theorem (\Cref{thm:hodge-thm}) identifies with the suffix points of $L^W$.}
  \label{fig:bridge}
\end{figure}


\subsection{Related Work}\label{sec:related}

Discrete Hodge theory~\citep{eckmann1944}, spectral graph theory~\citep{chung_spectral_1997}, and their extensions to discrete vector bundles~\citep{singer2012} and cellular sheaves~\citep{hansen2019,hansen2020} model consensus and signal processing on networks~\citep{olfati2007consensus,shuman2013emerging} via the graph, graph connection, and sheaf Laplacians.
These operators act on Hilbert spaces, as does the connection Laplacian of the flat vector bundle associated with a principal bundle~\citep{gao_geometry_2021}.
They rely on inner products and spectral decomposition, and so they do not apply directly to data that are logical, ordered, or metric.
We ground diffusion in enriched category theory instead, following the observation of \citet{Law73} that metric spaces are categories enriched in the Lawvere quantale.
Aggregation then becomes a weighted limit~\citep{Kelly2005}. The max-plus combinations of discrete event systems~\citep{butkovic2010} and the max-min combinations of fuzzy logic~\citep{ZADEH} are weighted joins in suitable $\quant$-categories.
The closest antecedent is the Tarski Laplacian of \citet{Ghrist_2022}, a local operator for cellular sheaves of complete lattices and join-preserving maps, which we recover over the Boolean quantale (\Cref{rem:tarski}).

In classical sheaf theory~\citep{bredon_sheaf_2012}, an assignment of data either is a global section or is not.
Approximate notions have appeared more recently.
The sheaf Laplacian of \citet{hansen2019} induces a quadratic form on the cochains that measures their disagreement across edges, and diffusion and smoothing minimize this form.
For sheaves of pseudometric spaces, \citet{Robinson2017} introduced the \emph{consistency radius}, the supremum over inclusions $U' \subseteq U$ of open sets of the distance between the value of an assignment on $U'$ and the restriction to $U'$ of its value on $U$~\citep[Definition~7]{Robinson2020}. The consistency radius quantifies the worst-case local disagreement and obstructs extension to a global section.
Our weighted global sections are related: over the Lawvere quantale, our stalks are Lawvere metric spaces and our restrictions are nonexpansive maps.
Distances need be neither symmetric nor finite, the maps are more restricted than the continuous maps of \citet{Robinson2020}, and the tolerance is prescribed separately for each edge and each direction.

\citet{baez1997} proposed an analogy between inner products and hom-objects, which becomes concrete for $\quant$-categories, whose hom-objects take values in a quantale such as $[0,\infty]$.
Once the inner product is replaced by the hom-object and linear adjoints by adjoint $\quant$-functors (\Cref{fig:bridge}), the Lawvere Laplacian categorifies the sheaf Laplacian at the level of $\quant$-categories.
Just as solutions of the heat equation on a Hilbert space converge to the kernel of the Laplacian, the harmonic flow on a complete $\quant$-category converges, under a continuity hypothesis, to a suffix point of the Lawvere Laplacian (\Cref{prop:kleene}).
Descending Kleene iteration and a generalization (\Cref{thm:TFPT}) of the fixed-point theorem of \citet{Tarski1955} take the place of the spectral theorem.


\subsection{Contributions}\label{sec:contributions}

Throughout, $\quant$ is a unital quantale, in general neither commutative nor affine, and $G$ is a finite simple graph.
Stated informally, our main results are the following.

\begin{mainthm}[Weighted global sections]\label{mthm:sections}
    Let $\shf$ be a cellular sheaf of $\quant$-categories on $G$, and let $W$ be a weight function, which assigns to each ordered pair $(u,v)$ of adjacent vertices a tolerance $W(u,v) \in \quant$ with $W(u,v) \odot W(v,u) \preceq 1$.
    Then the $\quant$-category $\wsections = [\widetilde{W},\shf]$ of natural transformations from the weight diagram $\widetilde{W}$ to $\shf$ is isomorphic to the full subcategory of the $\quant$-category of $0$-cochains (the tuples of vertex data) on those that agree up to $W$ across every edge (\Cref{thm:global-sections}).
    If $\quant$ is commutative, then $\wsections$ is the limit of $\shf$ weighted by $\widetilde{W}$ (\Cref{prop:weighted-limit}).
\end{mainthm}

\begin{mainthm}[Hodge--Lawvere]\label{mthm:hodge}
    Let $\shf$ and $W$ be as in \Cref{mthm:sections}.
    If the restrictions of $\shf$ have right adjoints and its vertex stalks are complete, then the $0$-cochains lying below their image under the Lawvere Laplacian $L^W$ are the weighted global sections (\Cref{thm:hodge-thm}).
    Under a continuity hypothesis, the harmonic flow converges to the greatest weighted global section below its initial $0$-cochain (\Cref{prop:kleene}).
    If, moreover, $W$ is central, for instance when $\quant$ is commutative, then the weighted global sections form a complete coreflective subcategory of the $0$-cochains, and, under the same continuity hypothesis, the harmonic flow from a $0$-cochain computes its coreflection (\Cref{prop:gamma-limits}, \Cref{cor:flow-coreflection}).
\end{mainthm}

For central weight functions, the existence theory for these equilibria is the enriched Knaster--Tarski theorem (\Cref{thm:TFPT}), by which the suffix, prefix, and fixed points of a $\quant$-endofunctor of a complete $\quant$-category form complete $\quant$-categories.
Over a noncommutative base, every weighted construction must fix a side of residuation, and the composition law of $\quant$-categories dictates it (\Cref{def:weighted}).
With these sides, the theorems above hold for all weight functions, except the enriched statements of \Cref{mthm:hodge}, which require central weight functions.
A synchronization problem over the quantale of relations exhibits a noncentral weight function for which the Lawvere Laplacian is monotone and the harmonic flow converges, although the Laplacian is not a $\quant$-functor (\Cref{ex:matchings}).
The harmonic flow specializes to the Bellman--Ford iteration (\Cref{ex:shortest-paths}) and to a path problem on the nerve of a covering that reconciles local cost estimates (\Cref{ex:map-merging}).
Quantale synchronization asks for values at the vertices that realize elements of $\quant$ prescribed along oriented edges (\Cref{ex:sync}).
An adjoint bisheaf with stalks $\Quant$ encodes this problem, its solutions up to a tolerance are the weighted global sections, and, under the continuity hypothesis, the harmonic flow computes the greatest solution below an initial estimate, as in clock synchronization (\Cref{ex:clocks}).
Max-plus discrete event systems yield bisheaves whose stalks are products of $\Quant^{\op}$ indexed by events (\Cref{ex:maxplus}).
Part of this work appears in the doctoral thesis of \citet{lopezEnrichedCategoriesQuantales2025}.


\subsection{Organization}\label{sec:organization}

We distinguish the levels of enrichment summarized in \Cref{tab:enrichment-hierarchy}, from truth values to the category of $\quant$-categories.
\Cref{sec:prelim} collects background on quantales, $\quant$-categories, weighted meets and joins, and fixed points.
\Cref{sec:sheaves} introduces cellular sheaves of $\quant$-categories and their weighted global sections.
\Cref{sec:lawvere-laplacian} constructs the Lawvere Laplacian, proves the Hodge--Lawvere theorem, and studies the harmonic flow and the coreflection onto the weighted global sections (\Cref{sec:structure}).
\Cref{sec:examples} instantiates the theory in shortest paths, the reconciliation of local cost estimates, quantale synchronization over commutative and noncommutative bases, and discrete event systems.
\Cref{sec:conclusion} records open directions: other domains, bases that vary from cell to cell, and the relation to linear Laplacians.
Appendix~\ref{appx:technical} collects the residuation calculus and lemmas on weighted meets and joins, and Appendix~\ref{appx:weighted-limit} establishes the universal property of the weighted global sections.

\begin{table}[htb]
    \centering
    \begin{tikzpicture}[
        x=1cm, y=1cm,
        head/.style={align=center, font=\scriptsize\bfseries, text=black, inner sep=1pt},
        symc/.style={align=center, font=\footnotesize, text=black, inner sep=1pt},
        cell/.style={align=center, font=\footnotesize, text=black, inner sep=1pt},
        rline/.style={black!28, line width=0.4pt},
      ]
      \useasboundingbox (0.00,-0.25) rectangle (16.00,3.78);

      \begin{scope}
        \clip[rounded corners=3pt] (0.00,0.00) rectangle (16.00,3.78);

        \fill[black!7]  (0.00,0.00)  rectangle (2.00,3.10);
        \fill[black!2]  (2.00,0.00)  rectangle (5.10,3.10);
        \fill[black!9]  (5.10,0.00)  rectangle (7.70,3.10);
        \fill[black!3]  (7.70,0.00)  rectangle (12.60,3.10);
        \fill[black!7]  (12.60,0.00) rectangle (16.00,3.10);

        \fill[black!15] (0.00,3.10)  rectangle (2.00,3.78);
        \fill[black!10] (2.00,3.10)  rectangle (5.10,3.78);
        \fill[black!18] (5.10,3.10)  rectangle (7.70,3.78);
        \fill[black!12] (7.70,3.10)  rectangle (12.60,3.78);
        \fill[black!15] (12.60,3.10) rectangle (16.00,3.78);

        \foreach \xx in {2.00, 5.10, 7.70, 12.60} {
          \draw[black!22, line width=0.35pt] (\xx,0.00) -- (\xx,3.78);
        }

        \draw[black!60, line width=0.65pt] (0.00,3.10) -- (16.00,3.10);
        \foreach \yy in {2.48, 1.86, 1.24, 0.62} {
          \draw[rline] (0.00,\yy) -- (16.00,\yy);
        }
      \end{scope}

      \draw[rounded corners=3pt, black!55, line width=0.55pt]
        (0.00,0.00) rectangle (16.00,3.78);

      \node[head] at (1.00,3.44)  {SYMBOL};
      \node[head] at (3.55,3.44)  {OBJECTS};
      \node[head] at (6.40,3.44)  {ENRICHED IN};
      \node[head] at (10.15,3.44) {HOM-OBJECTS};
      \node[head] at (14.30,3.44) {MONOIDAL PRODUCT};

      \node[symc] at (1.00,2.79)  {$\bool$};
      \node[cell] at (3.55,2.79)  {truth values};
      \node[cell] at (6.40,2.79)  {$\bool$};
      \node[cell] at (10.15,2.79) {$a \to b$};
      \node[cell] at (14.30,2.79) {$a \meet b$};

      \node[symc] at (1.00,2.17)  {$\quant$};
      \node[cell] at (3.55,2.17)  {elements of $\quant$};
      \node[cell] at (6.40,2.17)  {$\bool$};
      \node[cell] at (10.15,2.17) {$p \preceq q$};
      \node[cell] at (14.30,2.17) {$p \odot q$};

      \node[symc] at (1.00,1.55)  {$\Quant$};
      \node[cell] at (3.55,1.55)  {elements of $\quant$};
      \node[cell] at (6.40,1.55)  {$\quant$};
      \node[cell] at (10.15,1.55) {$[p,q]$};
      \node[cell] at (14.30,1.55) {$p \odot q$};

      \node[symc] at (1.00,0.93)  {$\C$};
      \node[cell] at (3.55,0.93)  {arbitrary};
      \node[cell] at (6.40,0.93)  {$\quant$};
      \node[cell] at (10.15,0.93) {$\hom_{\C}(x,y)$};

      \node[symc] at (1.00,0.31)  {$\QCat$};
      \node[cell] at (3.55,0.31)  {$\quant$-categories};
      \node[cell] at (6.40,0.31)  {preorders};
      \node[cell] at (10.15,0.31) {$\quant$-functors, forming $[\C,\D]$};
      \node[cell] at (14.30,0.31) {$\C \times \D$ (or $\C \tensor \D$)};
    \end{tikzpicture}
    \caption{Levels of enrichment, from truth values to the category of $\quant$-categories. Over a noncommutative base, $\Quant$ is enriched through the right hom $[p,q]$ (\Cref{ex:constructions}), and $\QCat$ is a locally preordered $2$-category that is enriched in $(\QCat,\tensor,\terminal)$ when $\quant$ is commutative (Appendix~\ref{appx:weighted-limit}).}
    \label{tab:enrichment-hierarchy}
\end{table}

\subsection*{Acknowledgements}

Riess was supported by the Defense Advanced Research Projects Agency (HR0011-25-3-0235) and the Air Force Office of Scientific Research (FA9550-26-1-B039).


\section{Preliminaries}\label{sec:prelim}

We fix notation for quantales, $\quant$-categories, weighted meets and joins, and fixed points.
Routine lemmas are collected in Appendix~\ref{appx:technical}.


\subsection{Quantales}\label{sec:quantales}

Coefficients are drawn from quantales.
We do not assume quantales to be commutative: the general theory is developed with the two residuals of the multiplication.

\begin{definition}\label{def:quantale}
    A \emph{unital quantale} is a complete lattice $(\quant, \preceq)$ endowed with an associative multiplication $\odot : \quant \times \quant \to \quant$, with a two-sided unit $1 \in \quant$, that distributes over arbitrary joins in each variable:
    \begin{equation*}
        q \odot \Bigl( \textstyle\bigjoin_{i \in I} q_i\Bigr) = \textstyle\bigjoin_{i \in I} (q \odot q_i),
        \qquad
        \Bigl( \textstyle\bigjoin_{i \in I} q_i\Bigr) \odot q = \textstyle\bigjoin_{i \in I} (q_i \odot q).
    \end{equation*}
    It is \emph{commutative} if $\odot$ is, and \emph{affine} if $1 = \top$.
\end{definition}

Since multiplication preserves joins in each variable, the maps $(-) \odot q$ and $q \odot (-)$ have right adjoints, the \emph{right hom} $[q,-]$ and the \emph{left hom} $\dlhom{q}{-}$, characterized by
\begin{equation}\label{eq:adjoint-hom}
    p \odot q \preceq r \iff p \preceq [q, r],
    \qquad\qquad
    q \odot p \preceq r \iff p \preceq \dlhom{q}{r}
\end{equation}
for all $p, q, r \in \quant$.
In words, $[q,r]$ is the largest $p$ with $p \odot q \preceq r$, and $\dlhom{q}{r}$ is the largest $p$ with $q \odot p \preceq r$.
When $\quant$ is commutative the two coincide, and we write $[q,r]$ for both.
\Cref{lem:functorality-of-hom} collects their calculus.

\begin{definition}\label{def:central}
    An element $c$ of a quantale $\quant$ is \emph{central} if $c \odot p = p \odot c$ for all $p \in \quant$.
\end{definition}

The elements $\bot$ and $1$ are central.
Central elements enter through \Cref{lem:lap-functor}, by which the Lawvere Laplacian is a $\quant$-functor whenever its weight function takes central values.

\begin{examples}\label{ex:quantales}
    \leavevmode
    \begin{itemize}
        \item[(a)] The \emph{Boolean quantale} $\bool = (\{\bot,\top\}, \meet, \top)$ is commutative and affine, and its hom is implication.
        \item[(b)] The \emph{Lawvere quantale} $\law = ([0,\infty], +, 0)$ of costs, ordered by $\geq$, is commutative and affine. Its joins are infima, and its hom is truncated subtraction, $[a, b] = \max(0, b-a)$, with $[\infty,\infty] = 0$.
        \item[(c)] For a set $X$, the \emph{quantale of relations} $\cat{Rel}(X)$ consists of the binary relations on $X$ ordered by inclusion, with relational composition $R \fatsemi S = \{(x,z) : (x,y) \in R \text{ and } (y,z) \in S \text{ for some } y \in X\}$ as multiplication and the diagonal $\Delta_X$ as unit. For $|X| \geq 2$ it is neither commutative nor affine. Its elements above the unit are the reflexive relations, and its central elements are $\emptyset$ and $\Delta_X$.
        \item[(d)] For a monoid $M$ with identity $e$, the powerset $\cat{P}(M)$, with multiplication $A \cdot B = \{ab : a \in A \text{ and } b \in B\}$ and unit $\{e\}$, is a unital quantale, commutative if and only if $M$ is. Formal languages under concatenation are the case of a free monoid.
        \item[(e)] The \emph{extended Lawvere quantale} $\elaw = ([-\infty,\infty], +, 0)$ is ordered by $\geq$, and its bottom element $\infty$ absorbs every element. It is commutative but not affine, since its top element $-\infty$ lies strictly above the unit $0$. Its hom is $[a,b] = b - a$ for finite $a$ and $b$, extended to $\pm\infty$ by residuation. Unlike the hom of $\law$, it is not truncated at $0$, and \Cref{ex:maxplus} relies on this absence of truncation.
    \end{itemize}
\end{examples}


\subsection{Quantale-Enriched Categories}\label{sec:q-categories}

\medskip\noindent We recall the basic notions of enriched category theory~\citep[Chapter~1]{Kelly2005} for a fixed unital quantale $\quant$.

\begin{definition}\label{def:q-category-main}
    A \emph{$\quant$-category} $\C$ consists of a set $\obj(\C)$ of objects and, for all $x, y \in \obj(\C)$, a \emph{hom-object} $\hom_{\C}(x,y) \in \quant$ such that
    \begin{itemize}
        \item[(a)] $\hom_{\C}(y,z) \odot \hom_{\C}(x,y) \preceq \hom_{\C}(x,z)$ for all $x,y,z \in \obj(\C)$, and
        \item[(b)] $1 \preceq \hom_{\C}(x,x)$ for all $x \in \obj(\C)$.
    \end{itemize}
\end{definition}

\medskip\noindent We write $x \in \C$ for $x \in \obj(\C)$.

\begin{notation}\label{not:qleq}
    For objects $x$ and $y$ of a $\quant$-category $\C$ and an element $q \in \quant$, we write $x \qleq_q y$ if $q \preceq \hom_\C(x,y)$, and $x \qequals_q y$ if $x \qleq_q y$ and $y \qleq_q x$.
    By \Cref{def:q-category-main}, the relation $\qleq_1$ is a preorder on $\obj(\C)$, the \emph{underlying preorder} of $\C$.
\end{notation}

\begin{definition}\label{def:q-functor}
    A \emph{$\quant$-functor} $f : \C \to \D$ is a map $f : \obj(\C) \to \obj(\D)$ such that $\hom_\C(x,y) \preceq \hom_\D(fx, fy)$ for all $x, y \in \C$.
    It is \emph{fully faithful} if equality holds for all $x,y \in \C$, and an \emph{isomorphism} if it is fully faithful and bijective on objects.
\end{definition}

A $\quant$-functor is monotone for the underlying preorders.
We write $\QCat$ for the category of small $\quant$-categories and $\quant$-functors.
The hom-objects of a $\quant$-category are elements of $\quant$, so that its underlying ordinary category is the preorder $\qleq_1$. Consequently, $\quant$-functors, and the Laplacians and flows of \Cref{sec:lawvere-laplacian}, are maps of objects, and the enrichment enters through inequalities between hom-objects.

\begin{examples}\label{ex:q-categories}
    \leavevmode
    \begin{itemize}
        \item[(a)] A $\bool$-category is a preorder, and a $\bool$-functor is a monotone map.
        \item[(b)] A $\law$-category is a \emph{Lawvere metric space}~\citep{Law73}, a set with a distance that vanishes on the diagonal and satisfies the triangle inequality but need be neither symmetric nor finite. A $\law$-functor is a nonexpansive map.
        \item[(c)] A $\cat{Rel}(X)$-category with object set $C$ amounts to a preorder $\sqsubseteq$ on $C \times X$, where $(a,b) \in \hom(x,y)$ if and only if $(y,a) \sqsubseteq (x,b)$. A $\cat{Rel}(X)$-functor $f$ is a map of objects such that $f \times \id_X$ is monotone.
        \item[(d)] If a monoid $M$ acts on a set $S$, then the transporters $\hom(x,y) = \{m \in M : m \cdot x = y\}$ make $S$ a $\cat{P}(M)$-category.
    \end{itemize}
\end{examples}

\begin{examples}\label{ex:constructions}
    \leavevmode
    \begin{itemize}
        \item[(a)] The quantale $\quant$ is a $\quant$-category $\Quant$ with $\hom_{\Quant}(p,q) = [p,q]$, the right hom. The composition law of \Cref{def:q-category-main}(a) is $[q,r] \odot [p,q] \preceq [p,r]$ (\Cref{lem:functorality-of-hom}(7)), and the underlying preorder is the order of $\quant$.
        \item[(b)] The \emph{product} $\prod_{i \in I} \C_i$ of $\quant$-categories has the tuples $(x_i)_{i \in I}$ as objects and the hom-objects
        \[ \hom\bigl((x_i)_{i \in I},(y_i)_{i \in I}\bigr) = \textstyle\bigmeet_{i \in I} \hom_{\C_i}(x_i,y_i). \]
        \item[(c)] For $\quant$-categories $\C$ and $\D$, the $\quant$-functors $\C \to \D$ form a $\quant$-category $[\C,\D]$ with
        \[ \hom_{[\C,\D]}(f,g) = \textstyle\bigmeet_{x \in \C} \hom_\D(fx,gx). \]
        \item[(d)] If $\quant$ is commutative, then the \emph{opposite} $\C^\op$ of a $\quant$-category $\C$ has the same objects and $\hom_{\C^\op}(x,y) = \hom_\C(y,x)$. Over a noncommutative base, $\C^\op$ is enriched in the reversed quantale, with multiplication $(p,q) \mapsto q \odot p$, and in general not in $\quant$.
        \item[(e)] The \emph{unit $\quant$-category} $\terminal$ has a single object $\ast$ with $\hom_{\terminal}(\ast,\ast) = 1$. A $\quant$-functor $\terminal \to \C$ amounts to an object of $\C$.
    \end{itemize}
\end{examples}


\subsection{Weighted Meets and Joins}\label{sec:completeness}

Classical Laplacians aggregate by sums, and quantale-enriched Laplacians aggregate by weighted meets and joins.
We need them only for families of objects indexed by a set.

\begin{definition}[{cf.~\citealp[Section~3.1]{Kelly2005}}]\label{def:weighted}
    Let $\C$ be a $\quant$-category, $\S$ a set, $D : \S \to \obj(\C)$ a family of objects, and $W : \S \to \quant$ a \emph{weight}.
    The \emph{meet of $D$ weighted by $W$} is an object $\bigmeet^W D$ of $\C$ such that
    \begin{equation*}
        \hom_{\C}\Bigl(x, \bigmeet\nolimits^W D\Bigr) = \bigmeet_{s \in \S} \dlhom{W(s)}{\hom_\C(x,D(s))} \qquad \text{for all } x \in \C,
    \end{equation*}
    and the \emph{join of $D$ weighted by $W$} is an object $\bigjoin^W D$ of $\C$ such that
    \begin{equation*}
        \hom_{\C}\Bigl(\bigjoin\nolimits^W D, x\Bigr) = \bigmeet_{s \in \S} \bigl[W(s), \hom_\C(D(s), x)\bigr] \qquad \text{for all } x \in \C.
    \end{equation*}
    A $\quant$-category is \emph{complete}, or a \emph{weighted lattice}, if it admits all weighted meets and joins, and a weighted meet or join is \emph{conical} if $W(s) = 1$ for all $s \in \S$.
\end{definition}

\medskip\noindent Weighted meets and joins are unique up to $\qequals_1$, and a conical meet or join is a meet or join in the underlying preorder.

The sides of the homs in \Cref{def:weighted} are dictated by the composition law of \Cref{def:q-category-main}(a), which makes every function of the form $\psi = \hom_\C(-,m)$ satisfy $\psi(y) \odot \hom_\C(x,y) \preceq \psi(x)$.
The right-hand side of the first formula satisfies this inequality for every weight, by the residuation~\eqref{eq:adjoint-hom} and the composition law, whereas the function $x \mapsto [W(s), \hom_\C(x,D(s))]$ does not in general.
For instance, for $\C = \Quant$ over $\cat{Rel}(\{0,1\})$, the latter function satisfies this inequality, and is then of the form $\hom_{\Quant}(-,m)$, for only $46$ of the $256$ pairs $\bigl(W(s), D(s)\bigr) \in \quant^2$, by exhaustive enumeration.
Dually, weighted joins use the right hom.
These notions are the weighted limits and colimits of \citet{stubbe_introduction_2014} for a one-object quantaloid (with composition read in diagrammatic order), and over a commutative base the two homs agree.

\begin{definition}[{cf.~\citealp[Section~3.7]{Kelly2005}}]\label{def:cotensor}
    A $\quant$-category $\C$ is \emph{cotensored} if for all $q \in \quant$ and $y \in \C$ there is an object $q \cotensor y$ of $\C$ such that $\hom_\C(x, q \cotensor y) = \dlhom{q}{\hom_\C(x,y)}$ for all $x \in \C$. It is \emph{tensored} if for all $q \in \quant$ and $x \in \C$ there is an object $q \tensor x$ of $\C$ such that $\hom_\C( q \tensor x, y) = [q, \hom_\C(x,y)]$ for all $y \in \C$.
\end{definition}

Cotensors and tensors are the weighted meets and joins of one-element families.
In a cotensored $\quant$-category with conical meets, and in a tensored $\quant$-category with conical joins, respectively, we have
\begin{equation}\label{eq:weighted-meet-formula}
    \bigmeet\nolimits^W D \qequals_1 \bigmeet_{s \in \S} W(s) \cotensor D(s),
    \qquad\qquad
    \bigjoin\nolimits^W D \qequals_1 \bigjoin_{s \in \S} W(s) \tensor D(s),
\end{equation}
since the two sides of each formula have the same hom-objects from or to every object~\citep[Section~3.10]{Kelly2005}.
The basic complete $\quant$-categories are built from $\quant$ itself.

\begin{lemma}\label{lem:self-enrichment}
    Let $\quant$ be a unital quantale.
    \begin{itemize}
        \item[(a)] The $\quant$-category $\Quant$ is complete. Its conical meets and joins are meets and joins in $\quant$, its cotensors are $q \cotensor y = \dlhom{q}{y}$, and its tensors are $q \tensor x = q \odot x$.
        \item[(b)] If $\quant$ is commutative, then $\Quant^{\op}$ is complete. Its conical meets and joins are joins and meets in $\quant$, its cotensors are $q \cotensor y = q \odot y$, and its tensors are $q \tensor x = [q, x]$.
        \item[(c)] A product of complete $\quant$-categories is complete, with weighted meets and joins computed componentwise.
    \end{itemize}
\end{lemma}

\medskip\noindent Completeness requires only half of \Cref{def:weighted}.

\begin{lemma}\label{lem:complete-cocomplete}
    A $\quant$-category that admits all weighted joins admits all weighted meets, and conversely.
\end{lemma}

\Cref{lem:complete-cocomplete} is classical~\citetext{\citealp[Proposition~5.10]{stubbe_categories_2005}; see also \citealp[Example~3.11]{stubbe_introduction_2014}}, and Appendix~\ref{appx:technical} proves both lemmas.
With the $\quant$-functors that preserve weighted joins (\Cref{sec:adjunctions}) as morphisms, complete $\quant$-categories correspond to the modules over $\quant$~\citep[Corollary~2.2]{stubbe2007modules}.


\subsection{Fixed Points}\label{sec:fixed-points}

The equilibria studied in the sequel are fixed points of maps, possibly one-sided.
Over $\quant = \bool$, \Cref{thm:TFPT} below is the Knaster--Tarski fixed-point theorem~\citep{knaster1928,Tarski1955}.

\begin{definition}\label{def:fixed-points}
    Let $\C$ be a $\quant$-category and $T : \obj(\C) \to \obj(\C)$ a map. The full subcategories $\cat{Prefix}(T)$, $\cat{Suffix}(T)$, and $\cat{Fix}(T)$ of $\C$ consist of the objects $x$ with $Tx \qleq_1 x$, with $x \qleq_1 Tx$, and with $x \qequals_1 Tx$, respectively. We call their objects the \emph{prefix points}, the \emph{suffix points}, and the \emph{fixed points} of $T$.
\end{definition}

\begin{lemma}\label{lem:limit-closure}
    Let $\C$ be a $\quant$-category and let $T : \C \to \C$ be a $\quant$-functor. Then $\cat{Prefix}(T)$ is closed in $\C$ under weighted meets, and $\cat{Suffix}(T)$ is closed in $\C$ under weighted joins.
\end{lemma}

\begin{proof}
    Let $D : \S \to \obj(\cat{Suffix}(T))$ and $W : \S \to \quant$, and suppose that $j = \bigjoin^W D$ exists in $\C$. For $s \in \S$ we have $W(s) \preceq \hom(D(s), j) \preceq \hom(TD(s), Tj)$, by \Cref{lem:weights} and functoriality, and $1 \preceq \hom(D(s), TD(s))$. The composition law gives
    \[ W(s) = W(s) \odot 1 \preceq \hom(T D(s), T j) \odot \hom(D(s), T D(s)) \preceq \hom(D(s), Tj), \]
    so $1 \preceq [W(s), \hom(D(s), Tj)]$ for every $s \in \S$, and $1 \preceq \hom(j, Tj)$ by \Cref{def:weighted}. For a weighted meet $m$ of prefix points, the same argument gives $W(s) \preceq \hom(TD(s), D(s)) \odot \hom(Tm, TD(s)) \preceq \hom(Tm, D(s))$, whence $1 \preceq \hom(Tm, m)$.
\end{proof}

\begin{lemma}\label{lem:endo-restricts}
    Let $\C$ be a $\quant$-category and let $T : \C \to \C$ be a $\quant$-functor. Then $T$ restricts to $\quant$-endofunctors of $\cat{Prefix}(T)$, $\cat{Suffix}(T)$, and $\cat{Fix}(T)$.
\end{lemma}

\begin{proof}
    By functoriality, $x \qleq_1 Tx$ implies $Tx \qleq_1 T^2x$, and $Tx \qleq_1 x$ implies $T^2 x \qleq_1 Tx$.
\end{proof}

\medskip\noindent \Cref{lem:complete-cocomplete,lem:limit-closure,lem:endo-restricts} combine into the enriched form of the Knaster--Tarski theorem.

\begin{theorem}[Enriched Knaster--Tarski theorem]\label{thm:TFPT}
    Let $\C$ be a complete $\quant$-category and $T : \C \to \C$ a $\quant$-functor. Then $\cat{Prefix}(T)$, $\cat{Suffix}(T)$, and $\cat{Fix}(T)$ are complete $\quant$-categories, and in particular nonempty.
\end{theorem}

\begin{proof}
    By \Cref{lem:limit-closure}, $\cat{Suffix}(T)$ admits all weighted joins, and hence all weighted meets by \Cref{lem:complete-cocomplete}. Dually, $\cat{Prefix}(T)$ is complete. By \Cref{lem:endo-restricts}, $T$ restricts to a $\quant$-endofunctor $T'$ of the complete $\quant$-category $\cat{Suffix}(T)$, and $\cat{Fix}(T) = \cat{Prefix}(T')$ is complete.
\end{proof}

\begin{remark}\label{rem:monotone-suffix}
    If $\C$ is complete and $T$ is only monotone for the underlying preorder, then the proof of \Cref{lem:limit-closure} still applies to conical joins, so $\cat{Suffix}(T)$ is closed under conical joins and has all meets in its underlying preorder.
\end{remark}


\section{Cellular Sheaves}\label{sec:sheaves}

Cellular sheaves on a graph assign data to its vertices and edges and maps to its incidences.
In this section we define them for $\quant$-categories and introduce their weighted global sections, the vertex data that agree across every edge up to prescribed tolerances.


\subsection{Sheaves on the Incidence Poset}\label{sec:incidence}

Throughout, $G$ is a finite simple graph with vertex set $V$ and edge set $E$, viewed as a $1$-dimensional regular cell complex whose $0$-cells are the vertices and whose $1$-cells are the edges.
We write $\partial e = \{v,w\}$ for the set of endpoints of an edge $e$, and the \emph{incidence poset} $\cat{G}$ of $G$ is the set of cells $V \sqcup E$ ordered by the face relation $v \face e$ whenever $v \in \partial e$.
For $v \in V$ we write $N(v)$ for the set of neighbors of $v$, and $w \sim_e v$ when $w \in N(v)$ and $\partial e = \{v,w\}$.
An index $w \sim_e v$ thus ranges over $N(v)$ and names the edge joining $w$ to $v$.

\begin{definition}\label{def:cellular-sheaf}
    A \emph{cellular sheaf of $\quant$-categories} on $G$, or \emph{cellular $\quant$-sheaf} for short, is a functor $\shf : \cat{G} \to \QCat$, and a \emph{cellular $\quant$-cosheaf} on $G$ is a functor $\cat{G}^{\op} \to \QCat$.
\end{definition}

Explicitly, a cellular $\quant$-sheaf $\shf$ assigns a $\quant$-category $\shf(\sigma)$, the \emph{stalk} at $\sigma$, to each cell $\sigma \in V \sqcup E$, and a $\quant$-functor $\shf_{v \face e} : \shf(v) \to \shf(e)$, the \emph{restriction} along $v \face e$, to each incidence $v \face e$.
No two incidences compose, so these data are subject to no further condition.
A cosheaf $\oshf$ assigns instead a $\quant$-functor $\oshf_{v \face e} : \oshf(e) \to \oshf(v)$ to each incidence $v \face e$. We index the maps of sheaves and of cosheaves by the same incidence $v \face e$, although they point in opposite directions.
\begin{remark}
    Our terminology follows the cellular-sheaf literature~\citep{curry2014,hansen2019}: the quantale $\quant$ is the base of enrichment of the stalks rather than a locale of open sets or its generalization, as in the theory of sheaves on a locale~\citep{walters1981} or on a quantale~\citep{borceux1986quantales}.
\end{remark}

\begin{remark}\label{rem:alexandrov}
    Cellular sheaves are sheaves in the topological sense~\citep{shepardCellularDescriptionDerived1985}.
    Endow $V \sqcup E$ with the Alexandrov topology, whose open sets are the up-sets of $\cat{G}$.
    Each cell $\sigma$ then has a smallest open neighborhood, its up-set or \emph{open star}, and sending $\sigma$ to its open star defines a functor from $\cat{G}$ to the opposite of the open-set lattice.
    Since $\QCat$ is complete and cocomplete as an ordinary category, right Kan extension along that functor gives an equivalence between functors $\shf : \cat{G} \to \QCat$ and $\QCat$-valued sheaves on $V \sqcup E$, taking $\shf$ to $U \mapsto \lim_{\sigma \in U} \shf(\sigma)$.
    For each cell $\sigma$, both the value of the induced sheaf on the open star of $\sigma$ and its topological stalk at $\sigma$ are canonically isomorphic to $\shf(\sigma)$.
    \citet[Theorem~4.2.10 and Corollary~4.2.13]{curry2014} stated the equivalence, and \citet[Theorem~4.8]{curry2019erratum} corrected its proof in the dual setting of cosheaves.
\end{remark}


\subsection{Weighted Global Sections}\label{sec:global-sections}

\medskip\noindent Throughout \Cref{sec:global-sections}, $\shf$ is a cellular $\quant$-sheaf on $G$.

\begin{definition}\label{def:cochains}
    The $\quant$-category of \emph{$0$-cochains} of $\shf$ is the product $\cochains \coloneqq \prod_{v \in V} \shf(v)$ of the vertex stalks. Its objects are the tuples $x = (x_v)_{v \in V}$, and its hom-objects are
    \[ \hom_{\cochains}(x, y) = \bigmeet_{v \in V} \hom_{\shf(v)}(x_v, y_v). \]
\end{definition}

We follow \citet[Section~3]{Ghrist_2022} in calling these $0$-cochains, although no coboundary operator appears below.
In the terminology of \citet{Robinson2020}, a $0$-cochain is an assignment supported on the open stars of the vertices.

\medskip\noindent Agreement across an edge is measured by a tolerance attached to each direction of the edge.

\begin{definition}\label{def:weight-function}
    A \emph{weight function} $W$ on $G$ assigns to each ordered pair $(u,v)$ of adjacent vertices an element $W(u,v) \in \quant$, the \emph{tolerance} from $u$ to $v$, such that $W(u,v) \odot W(v,u) \preceq 1$.
    It is \emph{symmetric} if $W(u,v) = W(v,u)$ for all adjacent vertices $u$ and $v$, and \emph{central} if all its values are central.
    For $q \in \quant$ with $q \odot q \preceq 1$, we write $W \equiv q$ for the weight function with constant value $q$.
\end{definition}

The larger the tolerance $W(u,v)$, the more closely the data at $u$ and $v$ must agree (\Cref{thm:global-sections}).
The tolerance $1$ in both directions across an edge demands agreement in the underlying preorder, and the tolerance $\bot$ imposes no constraint, so that a weight function with $W(u,v) = \bot$ and $W(v,u) \neq \bot$ imposes a one-sided condition across the edge.
The inequality $W(u,v) \odot W(v,u) \preceq 1$ is the composition law of \Cref{def:q-category-main}(a) for the two-object $\quant$-categories $\widetilde{W}(e)$ of \Cref{def:w-global-sections}, and it holds automatically when $\quant$ is affine.
Over a nonaffine base it restricts the tolerances that do not lie below the unit, and in particular those that lie above it.

\begin{lemma}\label{lem:above-unit}
    Let $W$ be a weight function on $G$ and let $u$ and $v$ be adjacent vertices. If $1 \preceq W(u,v)$, then $W(v,u) \preceq 1$. If moreover $1 \preceq W(v,u)$, then $W(u,v) = W(v,u) = 1$.
\end{lemma}

\begin{proof}
    Since $\odot$ is monotone, we have $W(v,u) = 1 \odot W(v,u) \preceq W(u,v) \odot W(v,u) \preceq 1$. If moreover $1 \preceq W(v,u)$, then $W(v,u) = 1$, and exchanging $u$ and $v$ gives $W(u,v) = 1$.
\end{proof}

Consequently, if $W(u,v)$ lies strictly above $1$, then $W(v,u)$ lies strictly below $1$.
Tolerances strictly above the unit, which demand more than agreement, thus exist only over a nonaffine base and only in one direction across an edge.

To describe the weighted global sections categorically, we encode a weight function as a diagram of $\quant$-categories and recall the natural transformations between such diagrams.
For functors $\Omega, \shf : \cat{G} \to \QCat$, a \emph{natural transformation} $\alpha : \Omega \Rightarrow \shf$ is a family of $\quant$-functors $\alpha_\sigma : \Omega(\sigma) \to \shf(\sigma)$ such that $\alpha_e \circ \Omega_{v \face e} = \shf_{v \face e} \circ \alpha_v$ for every incidence $v \face e$.
The natural transformations $\Omega \Rightarrow \shf$ form a $\quant$-category $[\Omega, \shf]$ with
\[ \hom_{[\Omega,\shf]}(\alpha, \beta) \coloneqq \bigmeet_{\sigma \in V \sqcup E} \hom_{[\Omega(\sigma), \shf(\sigma)]}(\alpha_\sigma, \beta_\sigma). \]

\begin{definition}\label{def:w-global-sections}
    Let $W$ be a weight function on $G$. The \emph{weight diagram} $\widetilde{W} : \cat{G} \to \QCat$ sends a vertex to the unit $\quant$-category $\terminal$, sends an edge $e$ with $\partial e = \{u,v\}$ to the $\quant$-category $\widetilde{W}(e)$ with objects $u$ and $v$ and hom-objects $\hom(u,v) = W(u,v)$, $\hom(v,u) = W(v,u)$, and $\hom(u,u) = \hom(v,v) = 1$, and sends an incidence $v \face e$ to the $\quant$-functor $\ast \mapsto v$. The $\quant$-category of \emph{$W$-weighted global sections} of $\shf$ is
    \[ \wsections \coloneqq [\widetilde{W}, \shf]. \]
\end{definition}

Over every unital quantale $\quant$, the $\quant$-category $\wsections = [\widetilde{W},\shf]$ is the end $\int_{\sigma \in \cat{G}} [\widetilde{W}(\sigma), \shf(\sigma)]$, constructed as an equalizer of products in the complete category $\QCat$. In Appendix~\ref{appx:weighted-limit}, we show that when $\quant$ is commutative, this end is the limit of $\shf$ weighted by $\widetilde{W}$ in the $\QCat$-enriched sense of \citet[Section~3.1]{Kelly2005}.
The following theorem describes $\wsections$ inside the $0$-cochains.

\begin{theorem}\label{thm:global-sections}
    Let $\shf$ be a cellular $\quant$-sheaf and $W$ a weight function on $G$. Evaluation at the vertices, $\alpha \mapsto (\alpha_v(\ast))_{v \in V}$, is an isomorphism of $\wsections$ onto the full subcategory of $\cochains$ on the $0$-cochains $x$ such that
    \begin{equation}\label{eq:w-section-condition}
        \shf_{v \face e} (x_v) \qleq_{W(v,w)} \shf_{w \face e} (x_w)
    \end{equation}
    for every edge $e$ and both orderings $(v,w)$ of its endpoints.
\end{theorem}

\begin{proof}
    A natural transformation $\alpha : \widetilde{W} \Rightarrow \shf$ consists of objects $x_v = \alpha_v(\ast)$ of the vertex stalks and $\quant$-functors $\alpha_e : \widetilde{W}(e) \to \shf(e)$ with $\alpha_e(v) = \shf_{v \face e}(x_v)$ for every $v \in \partial e$. The edge components are thus determined by the vertex components. For $\partial e = \{v,w\}$, the map $\alpha_e$ is a $\quant$-functor if and only if $W(v,w) \preceq \hom_{\shf(e)}(\shf_{v \face e}x_v, \shf_{w \face e}x_w)$ and $W(w,v) \preceq \hom_{\shf(e)}(\shf_{w \face e}x_w, \shf_{v \face e}x_v)$, the conditions at the identities of $\widetilde{W}(e)$ holding automatically. Evaluation is therefore a bijection onto the stated objects. For transformations $\alpha$ and $\beta$ with vertex components $x$ and $y$, the term of $\hom(\alpha,\beta)$ at an edge $e$ with $\partial e = \{v,w\}$ is $\hom(\shf_{v \face e}x_v, \shf_{v \face e}y_v) \meet \hom(\shf_{w \face e}x_w, \shf_{w \face e}y_w)$. It dominates $\hom(x_v,y_v) \meet \hom(x_w,y_w)$ because restrictions are $\quant$-functors, and hence $\hom(\alpha,\beta) = \bigmeet_{v \in V} \hom(x_v, y_v) = \hom_{\cochains}(x,y)$.
\end{proof}

We identify $\wsections$ with this full subcategory of $\cochains$.
The inequality $W(u,v) \odot W(v,u) \preceq 1$ of \Cref{def:weight-function} is used only to ensure that each $\widetilde{W}(e)$ is a $\quant$-category.
For an arbitrary assignment of tolerances $W(v,w) \in \quant$, the agreement condition~\eqref{eq:w-section-condition} still defines a full subcategory of $\cochains$. With this subcategory in place of $\wsections$, the Hodge--Lawvere theorem (\Cref{thm:hodge-thm}), the convergence of the harmonic flow (\Cref{prop:kleene}), and the coreflection (\Cref{prop:gamma-limits}) hold verbatim.
At $W \equiv 1$ the condition of \Cref{thm:global-sections} reads $\shf_{v \face e}x_v \qequals_1 \shf_{w \face e}x_w$, whereas the ordinary limit $\lim_{\cat{G}} \shf$, the global sections of the sheaf in \Cref{rem:alexandrov}, requires equality $\shf_{v \face e}x_v = \shf_{w \face e}x_w$ in each edge stalk.
The two coincide whenever the underlying preorders of the edge stalks are antisymmetric.


\section{The Lawvere Laplacian}\label{sec:lawvere-laplacian}

\Cref{thm:global-sections} characterizes the weighted global sections but provides no procedure that computes them.
In this section we construct an operator, the Lawvere Laplacian, whose suffix points are the weighted global sections (\Cref{thm:hodge-thm}), and a flow, the harmonic flow $x \mapsto x \meet L^W x$, that converges to a weighted global section under a continuity hypothesis (\Cref{prop:kleene}).


\subsection{Bisheaves and Transport}\label{sec:bisheaves}

\medskip\noindent The Lawvere Laplacian carries the datum at a vertex into an incident edge stalk and back to the other endpoint, so it requires maps in both directions along each incidence.

\begin{definition}\label{def:bisheaf}
    A \emph{cellular bisheaf} $\shf = (\ushf,\oshf)$ on $G$ consists of a cellular $\quant$-sheaf $\ushf$ and a cellular $\quant$-cosheaf $\oshf$ on $G$ with the same stalks, $\ushf(\sigma) = \oshf(\sigma)$ for all $\sigma \in V \sqcup E$.
\end{definition}

The restrictions of $\shf$ are those of $\ushf$, and its \emph{coextensions} are the maps $\oshf_{v \face e} : \shf(e) \to \shf(v)$.
The stalks $\shf(\sigma)$, the $0$-cochains $\cochains$, and the weighted global sections $\wsections$ of $\shf$ are those of its sheaf $\ushf$.
A restriction followed by a coextension transports data across an edge.

\begin{definition}\label{def:transport}
    Let $\shf$ be a bisheaf and $e$ an edge with $\partial e = \{v,w\}$. The \emph{transport} along $e$ from $w$ to $v$ is the $\quant$-functor
    \[ \transport{w}{v} \coloneqq \oshf_{v \face e} \circ \ushf_{w \face e} : \shf(w) \to \shf(v). \]
\end{definition}

\medskip\noindent The Lawvere Laplacian transports the data of the neighbors of a vertex into its stalk and aggregates them by a weighted meet (\Cref{fig:transport}).

\begin{definition}\label{def:lap}
    Let $\shf$ be a bisheaf and $W$ a weight function on $G$. We say that $\shf$ \emph{admits $W$-Laplacian meets} if for every vertex $v \in V$ and every $0$-cochain $x \in \obj(\cochains)$, the vertex stalk $\shf(v)$ admits the weighted meet
    \[ (L^W x)_v \coloneqq \bigmeet_{w \sim_e v}^{W(v,-)} \transport{w}{v}(x_w), \]
    the meet of the transported data weighted by $w \mapsto W(v,w)$, as holds whenever the vertex stalks are complete. The resulting map $L^W : \obj(\cochains) \to \obj(\cochains)$ is the \emph{Lawvere Laplacian} of $\shf$ weighted by $W$.
\end{definition}

Weighted meets are unique up to $\qequals_1$, and so is $L^W$. We fix a choice of representatives, on which none of the statements below depends. When the vertex stalks are cotensored, the decomposition~\eqref{eq:weighted-meet-formula} gives $(L^W x)_v \qequals_1 \bigmeet_{w \sim_e v} W(v,w) \cotensor \transport{w}{v}(x_w)$. At an isolated vertex $v$, the meet in \Cref{def:lap} is indexed by the empty set and equals the top object of $\shf(v)$, imposing no constraint on $x_v$.

\begin{figure}[tp]
  \centering
  \begin{tikzpicture}[
      x=1cm, y=1cm,
      panel/.style={rounded corners=3pt, line width=0.5pt},
      head/.style={inner sep=0pt, font=\footnotesize\bfseries, text=black},
      vtx/.style={circle, draw=black, fill=white, line width=0.65pt,
                  minimum size=6.2mm, inner sep=0pt, font=\footnotesize},
      centre/.style={vtx, fill=black!18, draw=black, line width=0.95pt},
      stalk/.style={rounded corners=3pt, draw=black!80, fill=white,
                    line width=0.65pt, minimum width=1.18cm,
                    minimum height=0.74cm, inner sep=2pt, font=\footnotesize},
      estalk/.style={rounded corners=2pt, draw=black!75, fill=black!10,
                     line width=0.55pt, inner sep=2.5pt, font=\scriptsize},
      map/.style={-{Straight Barb[length=3.4pt,width=3.8pt]},
                  black!80, line width=0.8pt},
      tmap/.style={-{Straight Barb[length=3.6pt,width=4.0pt]},
                   black, densely dashed, line width=0.85pt},
      lab/.style={font=\scriptsize, inner sep=1pt},
      wlab/.style={font=\scriptsize, inner sep=1pt, text=black!70},
      dlab/.style={font=\footnotesize, text=black, inner sep=1pt},
    ]

    \fill[black!4, panel, draw=black!30] (0.00,3.20) rectangle (7.75,9.05);
    \node[head, anchor=north west] at (0.28,8.80)
      {1.\enspace EDGE TRANSPORT AND ADJUNCTION};
    \draw[black!35, line width=0.4pt] (0.28,8.38) -- (7.47,8.38);

    \node[stalk] (FW) at (1.05,6.92) {$\shf(w)$};
    \node[stalk] (FE) at (3.88,6.92) {$\shf(e)$};
    \node[stalk] (FV) at (6.70,6.92) {$\shf(v)$};

    \draw[map] (FW.east) -- (FE.west)
      node[midway, above=2.5pt, font=\footnotesize, text=black] {$\ushf_{w \face e}$}
      node[midway, below=2.5pt, lab, text=black!65] {(restriction)};

    \draw[map] ([yshift=6.5pt]FV.west) -- ([yshift=6.5pt]FE.east)
      node[midway, above=2pt, font=\footnotesize, text=black] {$\ushf_{v \face e}$};
    \draw[map] ([yshift=-6.5pt]FE.east) -- ([yshift=-6.5pt]FV.west)
      node[midway, below=2pt, font=\scriptsize, align=center, inner sep=0pt]
        {{\footnotesize\textcolor{black}{$\oshf_{v \face e}$}}\\[0.5pt]
         \textcolor{black!65}{(coextension)}};
    \node[font=\scriptsize, text=black, inner sep=0pt]
      at ($(FE.east)!0.5!(FV.west)$)
      {$\rotatebox[origin=c]{-90}{$\dashv$}_{\mkern-1mu 1}$};

    \draw[tmap] ([xshift=-2pt]FW.south) to[out=-62, in=-118, looseness=1.14]
      node[pos=0.5, below=3pt, font=\footnotesize, text=black]
        {Transport\enspace $\transport{w}{v} = \oshf_{v \face e} \circ \ushf_{w \face e}$}
      ([xshift=2pt]FV.south);

    \draw[black!22, line width=0.35pt] (0.55,4.02) -- (7.20,4.02);
    \node[font=\scriptsize, text=black!75] at (3.88,3.58)
      {$\ushf_{v \face e} x_v \qleq_{W(v,w)} \ushf_{w \face e} x_w$
       \ $\Longleftrightarrow$ \
       $x_v \qleq_{W(v,w)} \transport{w}{v}(x_w)$};

    \fill[black!4, panel, draw=black!30] (8.25,3.20) rectangle (16.00,9.05);
    \node[head, anchor=north west] at (8.53,8.80)
      {2.\enspace LOCAL AGGREGATION AT A VERTEX};
    \draw[black!35, line width=0.4pt] (8.53,8.38) -- (15.72,8.38);

    \coordinate (vc) at (12.125,5.38);
    \coordinate (w1) at (12.125,7.82);
    \coordinate (w2) at (9.28,4.02);
    \coordinate (w3) at (14.97,4.02);

    \node[centre] (V)  at (vc) {$v$};
    \node[vtx]    (W1) at (w1) {$w_1$};
    \node[vtx]    (W2) at (w2) {$w_2$};
    \node[vtx]    (W3) at (w3) {$w_3$};

    \node[estalk] (E1) at ($(vc)!0.55!(w1)$) {$\shf(e_1)$};
    \node[estalk] (E2) at ($(vc)!0.54!(w2)$) {$\shf(e_2)$};
    \node[estalk] (E3) at ($(vc)!0.54!(w3)$) {$\shf(e_3)$};

    \draw[map] (W1) -- (E1);
    \draw[map] (E1) -- (V)
      node[wlab, midway, left=4pt]  {$W(v,w_1)$}
      node[lab,  midway, right=4pt, text=black] {$\transport{w_1}{v}(x_{w_1})$};

    \draw[map] (W2) -- (E2);
    \draw[map] (E2) -- (V)
      node[wlab, midway, above left=1.5pt and 1pt] {$W(v,w_2)$};
    \node[lab, text=black, anchor=north] at (10.95,4.22)
      {$\transport{w_2}{v}(x_{w_2})$};

    \draw[map] (W3) -- (E3);
    \draw[map] (E3) -- (V)
      node[wlab, midway, above right=1.5pt and 1pt] {$W(v,w_3)$};
    \node[lab, text=black, anchor=north] at (13.30,4.22)
      {$\transport{w_3}{v}(x_{w_3})$};

    \node[dlab, right=3.5pt of W1] {$x_{w_1}$};
    \node[dlab, below=2.5pt of W2] {$x_{w_2}$};
    \node[dlab, below=2.5pt of W3] {$x_{w_3}$};
    \node[dlab, below=3.5pt of V]  {$(L^W x)_v$};

    \fill[black!8, panel, draw=black!42] (1.35,-0.18) rectangle (14.65,2.78);
    \node[head, anchor=north west] at (1.65,2.53)
      {3.\enspace OPERATOR, FLOW, AND EQUILIBRIA};
    \draw[black!42, line width=0.4pt] (1.65,2.14) -- (14.35,2.14);

    \node[font=\footnotesize, text=black!70, anchor=base east] at (6.55,1.30)
      {Lawvere Laplacian:};
    \node[font=\footnotesize, text=black, anchor=base west]  at (6.80,1.30)
      {$\displaystyle (L^W x)_v = \bigmeet_{w \sim_e v}^{W(v,-)} \transport{w}{v}(x_w)$};
    \node[lab, text=black!70, anchor=base east] at (14.35,1.30) {\Cref{def:lap}};

    \node[font=\footnotesize, text=black!70, anchor=base east] at (6.55,0.60)
      {Harmonic flow:};
    \node[font=\footnotesize, text=black, anchor=base west]  at (6.80,0.60)
      {$x[t+1]_v = x[t]_v \meet (L^W x[t])_v$};
    \node[lab, text=black!70, anchor=base east] at (14.35,0.60) {\Cref{def:flow-functor}};

    \node[font=\footnotesize, text=black!70, anchor=base east] at (6.55,0.08)
      {Hodge--Lawvere theorem:};
    \node[font=\footnotesize, text=black, anchor=base west]  at (6.80,0.08)
      {$\cat{Suffix}(L^W) = \wsections$};
    \node[lab, text=black!70, anchor=base east] at (14.35,0.08) {\Cref{thm:hodge-thm}};

  \end{tikzpicture}
  \caption{Transport, local aggregation, and harmonic flow on a cellular bisheaf of $\quant$-categories. (1) On an edge $e$ with endpoints $v$ and $w$, the transport $\transport{w}{v}$ composes the restriction $\ushf_{w \face e}$ with the coextension $\oshf_{v \face e}$ (\Cref{def:transport}). When $\ushf_{v \face e} \dashv_1 \oshf_{v \face e}$, the adjunction transposes the agreement condition of \Cref{thm:global-sections} in $\shf(e)$ into an inequality in $\shf(v)$. (2) At a vertex $v$, the Lawvere Laplacian restricts the datum $x_{w_k}$ of each neighbor into $\shf(e_k)$, coextends it into $\shf(v)$, and forms the meet of the transported data weighted by the tolerances $W(v,w_k)$ (\Cref{def:lap}). (3) Iterating the binary meet with $L^W$ defines the harmonic flow (\Cref{def:flow-functor}), whose equilibria are the suffix points of $L^W$. For an adjoint bisheaf, the Hodge--Lawvere theorem (\Cref{thm:hodge-thm}) identifies them with the weighted global sections.}
  \label{fig:transport}
\end{figure}

\begin{lemma}\label{lem:lap-functor}
    Let $\shf$ be a bisheaf and $W$ a weight function on $G$ such that $\shf$ admits $W$-Laplacian meets. Then $L^W$ is monotone for the underlying preorder of $\cochains$. If, moreover, $W$ is central, for instance when $\quant$ is commutative, then $L^W$ is a $\quant$-functor.
\end{lemma}

\begin{proof}
    Let $x$ and $y$ be $0$-cochains and let $v \in V$ be a vertex. Transports are $\quant$-functors, so we have
    \[ \hom(x_w,y_w) \preceq \hom\bigl(\transport{w}{v}x_w, \transport{w}{v}y_w\bigr) \qquad \text{for every } w \sim_e v. \]
    If $x \qleq_1 y$, then $(L^W x)_v \qleq_1 (L^W y)_v$ by \Cref{lem:product-hom}(a). If $W$ is central, then \Cref{lem:product-hom}(b) gives $\hom(x,y) \preceq \hom((L^Wx)_v, (L^Wy)_v)$, and taking the meet over $v \in V$ gives $\hom(x,y) \preceq \hom(L^Wx, L^Wy)$.
\end{proof}

Centrality cannot be dropped. Over $\cat{Rel}(\{0,1\})$, on a single edge with stalks $\Quant$ and identity restrictions and coextensions, $L^W$ is a $\quant$-functor for only four of the $49$ weight functions, those with values in the center $\{\emptyset, \Delta_{\{0,1\}}\}$, by exhaustive enumeration.


\subsection{Adjunctions}\label{sec:adjunctions}

Hodge Laplacians are built from a coboundary operator and its adjoint.
The enriched analogue is an adjunction between restriction and coextension, which we also relax by an element of $\quant$.

\begin{definition}\label{def:fuzzy-adj}
    Let $q \in \quant$. A pair of $\quant$-functors $f : \C \rightleftarrows \D : g$ is \emph{$q$-adjoint}, written $f \dashv_q g$, if for all $x \in \C$ and $y \in \D$ we have
    \begin{equation}\label{eq:fuzzy-adjoint}
        \hom_\D(fx,y) \odot q \preceq \hom_\C(x,gy)
        \qquad\text{and}\qquad
        q \odot \hom_\C(x,gy) \preceq \hom_\D(fx,y).
    \end{equation}
\end{definition}

At $q = 1$ the two inequalities assert $\hom_\D(fx,y) = \hom_\C(x,gy)$, so a $1$-adjunction is a $\quant$-adjunction in the usual sense.
Over a commutative base, they read $\hom_\D(fx,y) \qequals_q \hom_\C(x,gy)$ in $\Quant$.
For $q \preceq 1$, a $q$-adjunction relaxes a $1$-adjunction, and each transposition of a hom-object multiplies it by $q$.

\begin{prop}\label{prop:adjoint-condition}
    Let $q \in \quant$. A pair of $\quant$-functors $f : \C \rightleftarrows \D : g$ is $q$-adjoint if and only if
    \[ x \qleq_q gf x \quad \text{and} \quad fgy \qleq_q y \qquad \text{for all } x \in \C \text{ and } y \in \D. \]
\end{prop}

\begin{proof}
    ($\Rightarrow$) Setting $y = fx$ in the first inequality of the $q$-adjunction condition~\eqref{eq:fuzzy-adjoint} gives $q = 1 \odot q \preceq \hom_\D(fx,fx) \odot q \preceq \hom_\C(x, gfx)$. Setting $x = gy$ in the second gives $q \preceq q \odot \hom_\C(gy,gy) \preceq \hom_\D(fgy, y)$.

    ($\Leftarrow$) By functoriality of $g$, the unit inequality, and the composition law, we have
    \[ \hom_\D(fx,y) \odot q \preceq \hom_\C(gfx, gy) \odot \hom_\C(x, gfx) \preceq \hom_\C(x, gy), \]
    and by functoriality of $f$ and the counit inequality, we have
    \[ q \odot \hom_\C(x,gy) \preceq \hom_\D(fgy, y) \odot \hom_\D(fx, fgy) \preceq \hom_\D(fx, y). \qedhere \]
\end{proof}

\begin{definition}\label{def:q-adjoint-bisheaf}
    Let $q \in \quant$. A bisheaf $\shf = (\ushf,\oshf)$ is \emph{$q$-adjoint} if $\ushf_{v \face e} \dashv_q \oshf_{v \face e}$ for every incidence $v \face e$, and \emph{adjoint} if it is $1$-adjoint.
\end{definition}

Let $\QAdj$ be the category whose objects are $\quant$-categories and whose morphisms $\C \to \D$ are $1$-adjunctions $f \dashv_1 g$ with $f : \C \to \D$, composed by $(f' \dashv_1 g') \circ (f \dashv_1 g) = (f'f \dashv_1 gg')$.
An adjoint bisheaf is then a functor $\cat{G} \to \QAdj$, which carries an incidence $v \face e$ to the adjunction $\ushf_{v \face e} \dashv_1 \oshf_{v \face e}$.

Adjoint bisheaves arise from a single $\quant$-functor per incidence: a restriction that preserves weighted joins has a right adjoint, which serves as the coextension.
A $\quant$-functor $f$ \emph{preserves weighted joins} if, for every family $D$ and weight $\varphi$ whose weighted join $\bigjoin^\varphi D$ exists, the object $f\bigl(\bigjoin^\varphi D\bigr)$ is a join of $fD$ weighted by $\varphi$, and preservation of weighted meets is defined dually.
Left adjoints preserve weighted joins, and right adjoints preserve weighted meets.

\begin{lemma}[{cf.~\citealp[Proposition~5.12]{stubbe_categories_2005}}]\label{lem:adjoints-preserve}
    Let $f : \C \rightleftarrows \D : g$ be $\quant$-functors with $f \dashv_1 g$. Then $g$ preserves weighted meets and $f$ preserves weighted joins.
\end{lemma}

\begin{proof}
    Let $D : \S \to \obj(\D)$ be a family and $\varphi : \S \to \quant$ a weight whose weighted meet exists in $\D$. For every $x \in \C$, the adjunction and \Cref{def:weighted} give
    \begin{align*}
        \hom_\C\Bigl(x, g \bigmeet\nolimits^\varphi D\Bigr) &= \hom_\D\Bigl(fx, \bigmeet\nolimits^\varphi D\Bigr) = \bigmeet_{s \in \S} \dlhom{\varphi(s)}{\hom_\D(fx, D(s))} \\
        &= \bigmeet_{s \in \S} \dlhom{\varphi(s)}{\hom_\C(x, gD(s))},
    \end{align*}
    so $g \bigmeet^\varphi D$ is a meet of $gD$ weighted by $\varphi$. The statement for $f$ is dual, with the right hom.
\end{proof}

\medskip\noindent The following form of the adjoint functor theorem is a converse, and it is used again in \Cref{sec:structure} and \Cref{ex:maxplus}.

\begin{prop}[{\citealp[Proposition~6.8]{stubbe_categories_2005}; cf.~\citealp[Theorem~2.4.3]{lopezEnrichedCategoriesQuantales2025}}]\label{prop:kan-bisheaf}
    Let $\C$ be a complete $\quant$-category, $\D$ a $\quant$-category, and $f : \C \to \D$ a $\quant$-functor that preserves weighted joins. Then $f$ has a right adjoint $f_{\ast} : \D \to \C$, given by the weighted join
    \[ f_{\ast} y \coloneqq \bigjoin\nolimits^{\hom_\D(f-,y)} \id_{\C} \]
    of the family $(x)_{x \in \C}$ weighted by $x \mapsto \hom_\D(fx, y)$.
\end{prop}

\begin{proof}
    The map $f_{\ast}$ is a $\quant$-functor: for all $y, y' \in \D$ and $x \in \C$, the composition law and \Cref{lem:weights} give $\hom_\D(y,y') \odot \hom_\D(fx,y) \preceq \hom_\D(fx,y') \preceq \hom_\C(x, f_{\ast}y')$, and taking the meet over $x$ gives $\hom_\D(y,y') \preceq \hom_\C(f_{\ast}y, f_{\ast}y')$ by \Cref{def:weighted}. We verify the conditions of \Cref{prop:adjoint-condition} at $q = 1$. By \Cref{lem:weights}, $1 \preceq \hom_\D(fx,fx) \preceq \hom_\C(x, f_{\ast}fx)$. Since $f$ preserves weighted joins, $f (f_{\ast} y) \qequals_1 \bigjoin^{\hom_\D(f-,y)} f$, and \Cref{def:weighted} gives
    \[ \hom_\D\Bigl(\bigjoin\nolimits^{\hom_\D(f-,y)} f, y\Bigr) = \bigmeet_{x \in \C} \bigl[\hom_\D(fx,y), \hom_\D(fx,y)\bigr] \succeq 1, \]
    so $f f_{\ast} y \qleq_1 y$. \Cref{prop:adjoint-condition} gives $f \dashv_1 f_{\ast}$.
\end{proof}

Dually, if $\D$ is a complete $\quant$-category, then every $\quant$-functor $g : \D \to \C$ that preserves weighted meets has a left adjoint $g_{!} : \C \to \D$, given by $g_{!} x \coloneqq \bigmeet^{\hom_\C(x,g-)} \id_{\D}$.
Consequently, restrictions from complete vertex stalks that preserve weighted joins, together with the coextensions $\oshf_{v \face e} \coloneqq (\ushf_{v \face e})_{\ast}$, form an adjoint bisheaf.

\medskip\noindent Adjoint functors often come in triples.

\begin{definition}\label{def:biadjoint}
    An adjoint bisheaf $\shf = (\ushf, \oshf)$ is \emph{biadjoint} if every restriction also has a left adjoint $\dshf_{v \face e} : \shf(e) \to \shf(v)$, so that every incidence carries an adjoint triple
    \[ \dshf_{v \face e} \dashv_1 \ushf_{v \face e} \dashv_1 \oshf_{v \face e}. \]
\end{definition}

\medskip\noindent By \Cref{prop:kan-bisheaf} and its dual, a $\quant$-functor between complete $\quant$-categories that preserves weighted meets and weighted joins has both adjoints.

\begin{example}[Bisheaves of presheaves]\label{ex:powers}
    Suppose that $\quant$ is commutative, and let $\sigma \mapsto \C_\sigma$ be a cellular $\quant$-cosheaf on $G$ with $\quant$-functors $g_{v \face e} : \C_e \to \C_v$.
    For each $\quant$-category $\C$, the presheaf $\quant$-category $[\C^{\op}, \Quant]$ (\Cref{ex:constructions}(c),(d)) is closed in $\prod_{\obj(\C)} \Quant$ under pointwise conical meets and joins, cotensors $[q,-]$, and tensors $q \odot (-)$ by \Cref{lem:functorality-of-hom}(5) and the distributivity of $\odot$, so it is complete by \Cref{lem:self-enrichment}(a),(c) and \eqref{eq:weighted-meet-formula}, with weighted meets and joins computed pointwise.
    Precomposition therefore defines restrictions $g_{v \face e}^{\ast} : [\C_v^{\op}, \Quant] \to [\C_e^{\op}, \Quant]$ that preserve weighted meets and joins, and by \Cref{prop:kan-bisheaf} and its dual their enriched Kan extensions $(g_{v \face e})_{!} \dashv_1 g_{v \face e}^{\ast} \dashv_1 (g_{v \face e})_{\ast}$ assemble into a biadjoint bisheaf with coextensions $(g_{v \face e})_{\ast}$.
    When each $\C_\sigma$ is the discrete $\quant$-category on a set $Y_\sigma$ (with $\hom(y,y) = 1$ and $\hom(y,y') = \bot$ for $y \neq y'$), the stalks reduce to the powers $\Quant^{Y_\sigma} = \prod_{Y_\sigma} \Quant$, and precomposition along maps of sets $g_{v \face e} : Y_e \to Y_v$ yields a biadjoint \emph{bisheaf of powers} over every unital quantale $\quant$, commutative or not.
    Over $\quant = \bool$, general preorders $\C_\sigma$ give stalks of down-sets, while discrete preorders give the \emph{powersheaves} with stalks $2^{Y_\sigma}$ and adjoint triples $(g_{v \face e})_{!} \dashv_1 g_{v \face e}^{-1} \dashv_1 (g_{v \face e})_{\ast}$ formed by direct image, preimage, and universal image.
\end{example}

\medskip\noindent Starting instead from maps $f_{v \face e} : \C_v \to \C_e$ gives adjunctions $(f_{v \face e})_{!} \dashv_1 f_{v \face e}^{\ast}$ on presheaves whose restrictions $(f_{v \face e})_{!}$ need not preserve meets (such as direct image along a noninjective map over $\bool$), and the synchronization and max-plus bisheaves of \Cref{ex:sync,ex:maxplus} are likewise adjoint without being biadjoint in general.


\subsection{The Hodge--Lawvere Theorem}\label{sec:hodge}

\medskip\noindent The following theorem is a degree-zero analogue of the cellular Hodge theorem ($\ker L \cong H^0(G;\shf)$, \Cref{fig:bridge}): it identifies the suffix points of the Lawvere Laplacian with the weighted global sections.

\begin{theorem}[Hodge--Lawvere theorem]\label{thm:hodge-thm}
    Let $\shf$ be an adjoint bisheaf and $W$ a weight function on $G$ such that $\shf$ admits $W$-Laplacian meets. Then $\cat{Suffix}(L^W) = \wsections$.
\end{theorem}

\begin{proof}
    Both sides are full subcategories of $\cochains$, so it suffices to compare objects. Let $x$ be a $0$-cochain and $v \in V$. \Cref{def:weighted} gives
    \begin{equation}\label{eq:hodge-unfold}
        \hom_{\shf(v)}\bigl(x_v, (L^W x)_v\bigr) = \bigmeet_{w \sim_e v} \dlhom{W(v,w)}{\hom_{\shf(v)}(x_v, \oshf_{v \face e} \ushf_{w \face e} x_w)},
    \end{equation}
    and $1 \preceq \dlhom{q}{p}$ if and only if $q \preceq p$ (\Cref{lem:functorality-of-hom}(4)). Hence $x \qleq_1 L^W x$ if and only if
    \[ W(v,w) \preceq \hom_{\shf(v)}(x_v, \oshf_{v \face e} \ushf_{w \face e} x_w) \qquad \text{for all } v \in V \text{ and } w \sim_e v. \]
    The adjunction $\ushf_{v \face e} \dashv_1 \oshf_{v \face e}$ identifies the right-hand side with $\hom_{\shf(e)}(\ushf_{v \face e} x_v, \ushf_{w \face e} x_w)$, and the condition becomes that of \Cref{thm:global-sections}.
\end{proof}

If the bisheaf is adjoint and its vertex stalks are complete, then the bottom $0$-cochain, the conical join of the empty family, lies below its image under $L^W$, and it is therefore a weighted global section for every weight function.
More generally, $\wsections = \cat{Suffix}(L^W)$ is then closed under conical joins and has all meets in its underlying preorder, by \Cref{rem:monotone-suffix}, since $L^W$ is monotone.

The proof of \Cref{thm:hodge-thm} transposes a single hom-object along an adjunction.
For a $q$-adjoint bisheaf, each transposition multiplies by $q$, and the identification weakens to two containments.
For a weight function $W$ on $G$ and an element $q \preceq 1$ of $\quant$, we write $W \odot q$ and $q \odot W$ for the functions with values $W(v,w) \odot q$ and $q \odot W(v,w)$, which are weight functions since $q \preceq 1$.

\begin{corollary}\label{cor:hodge-approx}
    Let $q \in \quant$ with $q \preceq 1$, let $\shf$ be a $q$-adjoint bisheaf, and let $W$ be a weight function on $G$ such that $\shf$ admits $W$-Laplacian meets and $(W \odot q)$-Laplacian meets. Then
    \[ \wsections \subseteq \cat{Suffix}(L^{W \odot q}) \qquad \text{and} \qquad \cat{Suffix}(L^W) \subseteq \cat{\Gamma}^{q \odot W}(G;\shf). \]
\end{corollary}

\begin{proof}
    Let $x$ be a $0$-cochain, and for $v \in V$ and $w \sim_e v$ let $a_{vw} = \hom_{\shf(v)}(x_v, \oshf_{v \face e}\ushf_{w \face e}x_w)$ and $b_{vw} = \hom_{\shf(e)}(\ushf_{v \face e}x_v, \ushf_{w \face e}x_w)$. The $q$-adjunction gives $b_{vw} \odot q \preceq a_{vw}$ and $q \odot a_{vw} \preceq b_{vw}$. For every weight function $W'$ for which $\shf$ admits $W'$-Laplacian meets, the unfolding~\eqref{eq:hodge-unfold} shows that $x \in \cat{Suffix}(L^{W'})$ if and only if $W'(v,w) \preceq a_{vw}$ for all $v$ and $w \sim_e v$, and by \Cref{thm:global-sections}, $x \in \cat{\Gamma}^{W'}(G;\shf)$ if and only if $W'(v,w) \preceq b_{vw}$ for all $v$ and $w \sim_e v$. If $W(v,w) \preceq b_{vw}$, then $W(v,w) \odot q \preceq b_{vw} \odot q \preceq a_{vw}$, and if $W(v,w) \preceq a_{vw}$, then $q \odot W(v,w) \preceq q \odot a_{vw} \preceq b_{vw}$.
\end{proof}

\medskip\noindent Over a commutative base both weight functions in \Cref{cor:hodge-approx} are $q \odot W$, and at $q = 1$ the corollary is \Cref{thm:hodge-thm}.

A $q$-adjunction need not be an adjunction. Over $\law$, for costs $\tau, \tau' \in [0,\infty)$ and $q \in [0,\infty]$, the $\law$-functors $(-) + \tau$ and $[\tau', -]$ on $\Quant$ are $q$-adjoint if and only if $q \geq |\tau - \tau'|$ in the numerical order (\Cref{prop:adjoint-condition}), and they are adjoint if and only if $\tau' = \tau$.

\begin{remark}[Boolean coefficients]\label{rem:tarski}
    Over $\quant = \bool$, a complete $\bool$-category is a complete lattice up to isomorphism of objects, an adjoint bisheaf equips each join-preserving restriction $\ushf_{v \face e}$ with its right adjoint $\oshf_{v \face e}$, and the weighted meet at $W \equiv \top$ is the meet. Writing $L^{\top}$ for $L^{W \equiv \top}$, the Lawvere Laplacian becomes
    \[ (L^{\top} x)_v = \bigmeet_{w \sim_e v} \oshf_{v \face e}\bigl(\ushf_{w \face e}(x_w)\bigr), \]
    the Tarski Laplacian of \citet[Definition~3.1]{Ghrist_2022}, whose suffix points are the fixed points of the harmonic flow $x \mapsto x \meet L^{\top} x$~\citep[Definition~3.4]{Ghrist_2022}. The meet in \citet{Ghrist_2022} also runs over $w = v$, but the extra term $\oshf_{v \face e}\ushf_{v \face e}(x_v)$ lies above $x_v$ by the unit of the adjunction, and so it changes neither the suffix points nor the flow. \Cref{thm:hodge-thm} specializes to the theorem of \citet[Theorem~3.1, Corollary~3.1]{Ghrist_2022} that identifies these fixed points with the global sections, and \Cref{thm:TFPT} specializes to the fixed-point theorem of \citet{Tarski1955}.
\end{remark}


\subsection{Harmonic Flow}\label{sec:flow}

\medskip\noindent Throughout \Cref{sec:flow}, $\shf$ is a bisheaf with complete vertex stalks and $W$ is a weight function on $G$.

\begin{definition}\label{def:flow-functor}
    Let $\Phi$ be the map on $0$-cochains with components
    \[ (\Phi x)_v \coloneqq x_v \meet (L^W x)_v, \]
    the binary conical meet of $x_v$ and $(L^W x)_v$. The \emph{harmonic flow} from a $0$-cochain $x[0]$ is the sequence $x[t+1] \coloneqq \Phi(x[t])$.
\end{definition}

By \Cref{lem:lap-functor}, $\Phi$ is monotone.
A $0$-cochain $x$ is an \emph{equilibrium} of the harmonic flow if $\Phi x \qequals_1 x$, which holds if and only if $x \qleq_1 L^W x$, so the equilibria are the suffix points of $L^W$.
The flow descends, $x[t+1] \qleq_1 x[t]$, and we write $x[\infty] \coloneqq \bigmeet_{t \geq 0} x[t]$ for the conical meet of its iterates, which the flow need not reach in finitely many steps.
Monotonicity of $L^W$ gives $L^W(x[\infty]) \qleq_1 \bigmeet_t L^W(x[t])$.
We call the reverse inequality the \emph{continuity hypothesis} for the flow from $x[0]$.
Under it, the next proposition shows by descending Kleene iteration~\citep{Cousot79} that $x[\infty]$ is an equilibrium.

\begin{prop}\label{prop:kleene}
    Let $\shf$ be a bisheaf with complete vertex stalks, $W$ a weight function on $G$, and $x[t]$ the harmonic flow from a $0$-cochain $x[0]$.
    \begin{itemize}
        \item[(a)] If $\bigmeet_{t} L^W(x[t]) \qleq_1 L^W(x[\infty])$, then $x[\infty] \in \cat{Suffix}(L^W)$.
        \item[(b)] If, in addition, $\shf$ is adjoint, then $x[\infty]$ is the greatest weighted global section below $x[0]$ in the underlying preorder.
    \end{itemize}
\end{prop}

\begin{proof}
    (a) For every $t$, we have $x[\infty] \qleq_1 x[t+1] \qleq_1 L^W x[t]$. Hence $x[\infty] \qleq_1 \bigmeet_t L^W(x[t]) \qleq_1 L^W(x[\infty])$.

    (b) By (a) and \Cref{thm:hodge-thm}, $x[\infty]$ is a weighted global section below $x[0]$. Let $y$ be a weighted global section with $y \qleq_1 x[0]$. If $y \qleq_1 x[t]$, then $y \qleq_1 L^W y \qleq_1 L^W x[t]$ by \Cref{thm:hodge-thm} and \Cref{lem:lap-functor}, so $y \qleq_1 x[t+1]$. By induction, $y \qleq_1 x[t]$ for all $t$, and hence $y \qleq_1 x[\infty]$.
\end{proof}

\begin{remark}\label{rem:chain-preservation}
    The hypothesis of \Cref{prop:kleene}(a) holds in each of the following cases.
    \begin{itemize}
        \item[(i)] The flow stabilizes, $x[t+1] \qequals_1 x[t]$ for some $t$, and then $x[\infty] \qequals_1 x[t]$. The flow stabilizes if every vertex stalk has finitely many objects up to $\qequals_1$.
        \item[(ii)] The bisheaf is biadjoint (\Cref{def:biadjoint}). Restrictions and coextensions are then right adjoints and preserve weighted meets (\Cref{lem:adjoints-preserve}), and in particular conical meets. The transports therefore preserve conical meets, and since weighted meets commute with conical meets (\Cref{lem:functorality-of-hom}(2)), $L^W$ preserves all conical meets.
        \item[(iii)] The Laplacian $L^W$ preserves the conical meets of descending sequences, as in \Cref{ex:maxplus}. Case~(ii) is an instance.
    \end{itemize}

    Without the hypothesis of \Cref{prop:kleene}(a), the limit $x[\infty]$ need not be a weighted global section. Over $\quant = \bool$ with $W \equiv \top$ on a single edge $e$ between $u$ and $v$, let $\mathbb{N} = \{0, 1, 2, \dots\}$ and let $\shf(u) = 2^{\mathbb{N}}$, $\shf(v) = 2^{\mathbb{N}} \times \bool$, and $\shf(e) = 2^{\mathbb{N}} \times 2^{\mathbb{N}} \times \bool$, ordered componentwise by inclusion and the order of $\bool$. Write $f(S) \coloneqq \{n + 1 : n \in S\}$ for the forward shift on $2^{\mathbb{N}}$, which preserves arbitrary unions and has right adjoint $f_{\ast}(B) = \{n \in \mathbb{N} : n + 1 \in B\}$, and write $\delta(S) \coloneqq \top$ if $S \neq \emptyset$ and $\delta(\emptyset) \coloneqq \bot$. The restrictions $\ushf_{u \face e}(S) \coloneqq \bigl(S, f(S), \delta(S)\bigr)$ and $\ushf_{v \face e}(T, b) \coloneqq \bigl(f(T), T, b\bigr)$ preserve arbitrary joins, with right adjoints $\oshf_{u \face e}(A, B, \top) = A \cap f_{\ast}(B)$, $\oshf_{u \face e}(A, B, \bot) = \emptyset$, and $\oshf_{v \face e}(A, B, b) = \bigl(f_{\ast}(A) \cap B, b\bigr)$. Starting from $x[0]_u = \mathbb{N}$ and $x[0]_v = (\mathbb{N}, \top)$, induction gives $x[t]_u = S_t$ and $x[t]_v = (S_t, \top)$ for every $t \geq 0$, where $S_t \coloneqq \{n \in \mathbb{N} : n \geq t\}$, since $f(S_t) \cap f_{\ast}(S_t) = S_{t+1}$ and $\delta(S_t) = \top$. Taking the meet over $t \geq 0$ yields $x[\infty]_u = \emptyset$ and $x[\infty]_v = (\emptyset, \top)$, which is not a global section because $\ushf_{u \face e}(\emptyset) = (\emptyset, \emptyset, \bot) \neq (\emptyset, \emptyset, \top) = \ushf_{v \face e}(\emptyset, \top)$.
\end{remark}

\medskip\noindent If $L^W$ preserves binary conical meets, as it does in case~(ii), then the iterates of the flow are meets of powers of the Laplacian.

\begin{corollary}\label{cor:kleene-powers}
    Let $\shf$ be a bisheaf with complete vertex stalks, $W$ a weight function on $G$, and $x[t]$ the harmonic flow from a $0$-cochain $x[0]$. If $L^W$ preserves binary conical meets, then $x[t] \qequals_1 \bigmeet_{k = 0}^{t} (L^W)^{k}(x[0])$ for every $t$, and consequently $x[\infty] \qequals_1 \bigmeet_{k \geq 0} (L^W)^{k}(x[0])$.
\end{corollary}

\begin{proof}
    Let $M[t] \coloneqq \bigmeet_{k \le t}(L^W)^k(x[0])$. If $x[t] \qequals_1 M[t]$, then $x[t+1] \qequals_1 M[t] \meet L^W(M[t])$, and preservation of binary meets gives $L^W(M[t]) \qequals_1 \bigmeet_{k \le t} (L^W)^{k+1}(x[0])$, whence $x[t+1] \qequals_1 M[t+1]$. Taking conical meets over $t$ gives the second claim.
\end{proof}


\subsection{Weighted Global Sections as a Coreflective Subcategory}\label{sec:structure}

The harmonic flow from a $0$-cochain converges, under a continuity hypothesis, to the greatest weighted global section below it (\Cref{prop:kleene}).
We next determine the weighted meets and joins of $\wsections$.
Its weighted joins are those of the $0$-cochains, and its weighted meets are obtained by applying a coreflector to those of the $0$-cochains.

\begin{prop}\label{prop:gamma-limits}
    Let $\shf$ be an adjoint bisheaf with complete vertex stalks, let $W$ be a central weight function on $G$, and let $i : \wsections \to \cochains$ be the inclusion.
    \begin{itemize}
        \item[(a)] The $\quant$-category $\wsections$ is complete, and it is closed in $\cochains$ under weighted joins, which are computed componentwise.
        \item[(b)] The inclusion $i$ has a right adjoint $\rho$, the \emph{coreflector}: for every $y \in \wsections$ and $x \in \cochains$, we have $\hom_{\cochains}(i(y), x) = \hom_{\wsections}(y, \rho(x))$, and $\rho(x)$ is the greatest weighted global section below $x$.
        \item[(c)] For every family $D : \S \to \obj(\wsections)$ and weight $\varphi : \S \to \quant$, the weighted meet of $D$ in $\wsections$ is $\rho\bigl(\bigmeet^\varphi iD\bigr)$, where $\bigmeet^\varphi iD$ is computed componentwise in $\cochains$.
        \item[(d)] If $\shf$ is biadjoint (\Cref{def:biadjoint}), then $\wsections$ is closed in $\cochains$ under weighted meets.
    \end{itemize}
\end{prop}

\begin{proof}
    (a) By \Cref{lem:lap-functor}, $L^W$ is a $\quant$-functor, so $\wsections = \cat{Suffix}(L^W)$ (\Cref{thm:hodge-thm}) is complete by \Cref{thm:TFPT} and closed under weighted joins by \Cref{lem:limit-closure}. Weighted joins in the product $\cochains$ are componentwise by \Cref{lem:self-enrichment}(c).

    (b) By (a), the inclusion $i$ preserves weighted joins, and \Cref{prop:kan-bisheaf} supplies a right adjoint $\rho$ with $i \dashv_1 \rho$. In the underlying preorders, the adjunction says that $y \qleq_1 \rho(x)$ if and only if $y \qleq_1 x$, for every weighted global section $y$.

    (c) The right adjoint $\rho$ preserves weighted meets (\Cref{lem:adjoints-preserve}). Since $i$ is fully faithful, $\rho i (y) \qequals_1 y$ for every weighted global section $y$, and applying $\rho$ to $\bigmeet^\varphi iD$ gives $\bigmeet^\varphi \rho i D \qequals_1 \bigmeet^\varphi D$.

    (d) Let $m \coloneqq \bigmeet^\varphi iD$. Each $D(s)$ is a suffix point of $L^W$, so $m \qleq_1 \bigmeet^\varphi L^W iD$ by \Cref{lem:product-hom}(a). For a biadjoint bisheaf, restrictions and coextensions preserve weighted meets (\Cref{lem:adjoints-preserve}). For $v \in V$ and $z \in \shf(v)$, let $h_{ws}(z) \coloneqq \hom(z, \transport{w}{v} D(s)_w)$. Since weighted meets of $0$-cochains are componentwise (\Cref{lem:self-enrichment}(c)) and the transports preserve them, \Cref{def:weighted} and \Cref{lem:functorality-of-hom}(2) and~(6) give
    \begin{align*}
        \hom\bigl(z, (L^W m)_v\bigr) &= \bigmeet_{w \sim_e v} \bigmeet_{s \in \S} \dlhom{\varphi(s) \odot W(v,w)}{h_{ws}(z)}, \\
        \hom\Bigl(z, \Bigl(\bigmeet\nolimits^\varphi L^W iD\Bigr)_v\Bigr) &= \bigmeet_{w \sim_e v} \bigmeet_{s \in \S} \dlhom{W(v,w) \odot \varphi(s)}{h_{ws}(z)},
    \end{align*}
    which agree since $W$ is central. Hence $m \qleq_1 L^W m$, and $m$ is a weighted global section.
\end{proof}

\medskip\noindent Under the hypothesis of \Cref{prop:kleene}(a), the harmonic flow from a $0$-cochain $x[0]$ converges to its coreflection $\rho(x[0])$.

\begin{corollary}\label{cor:flow-coreflection}
    Let $\shf$ and $W$ be as in \Cref{prop:gamma-limits}, and let $x[t]$ be the harmonic flow from a $0$-cochain $x[0]$. If $\bigmeet_{t} L^W(x[t]) \qleq_1 L^W(x[\infty])$, then $x[\infty] \qequals_1 \rho(x[0])$.
\end{corollary}

\begin{proof}
    By \Cref{prop:kleene}(b), $x[\infty]$ is the greatest weighted global section below $x[0]$ in the underlying preorder, and by \Cref{prop:gamma-limits}(b) so is $\rho(x[0])$.
\end{proof}

\begin{example}\label{ex:meet-not-pointwise}
    To see that biadjointness is needed in \Cref{prop:gamma-limits}(d), let $\quant = \bool$ and $W \equiv \top$ on a single edge $e$ between $u$ and $v$. Take $\shf(u) = 2^{\{1,2\}}$ ordered by inclusion and $\shf(v) = \shf(e) = \bool$, with $\ushf_{u \face e}(S) = \top$ if $S \neq \emptyset$ and $\bot$ if $S = \emptyset$ (whose right adjoint sends $\top$ to $\{1,2\}$ and $\bot$ to $\emptyset$), and $\ushf_{v \face e} = \oshf_{v \face e} = \id_{\bool}$. Then $\shf$ is an adjoint bisheaf, but $\ushf_{u \face e}$ does not preserve binary meets, so $\shf$ is not biadjoint. The global sections $(\{1\},\top)$ and $(\{2\},\top)$ have componentwise meet $(\emptyset,\top)$ in $\cochains$, which is not a global section, and their meet in $\wsections$ is $\rho(\emptyset,\top) = (\emptyset,\bot)$.
\end{example}


\subsection{Examples}\label{sec:examples}

The following examples instantiate the theory. They are shortest paths, a path problem on the nerve of a covering, quantale synchronization with its instances over clocks and over binary relations, and the synchronization of coupled max-plus systems.

\begin{example}[Shortest paths]\label{ex:shortest-paths}
    Let $\quant = \law$, and let every stalk be $\Quant^{\op}$. The objects of $\Quant^{\op}$ are the costs $x \in [0,\infty]$, its underlying order is the numerical order, its cotensors are $q \cotensor y = q + y$, and its conical meets are numerical minima (\Cref{lem:self-enrichment}(b)). Let every restriction and coextension be the identity. The resulting bisheaf is biadjoint. Assign to each edge $e$ a length $\ell(e) \in [0,\infty)$, and let $W(v,w) = W(w,v) = \ell(e)$ for $\partial e = \{v,w\}$. By \Cref{thm:global-sections}, a $0$-cochain $x$ is a weighted global section if and only if $x_v \leq x_w + \ell(e)$ for every edge $e$ and both orderings $(v,w)$ of its endpoints, that is, if $x$ is $1$-Lipschitz for the edge lengths. The Lawvere Laplacian
    \[ (L^W x)_v = \min_{w \sim_e v} \bigl(\ell(e) + x_w\bigr) \]
    is the Bellman operator of the shortest-path problem, and the harmonic flow $x[t+1]_v = \min\bigl(x[t]_v, (L^W x[t])_v\bigr)$ is the Bellman--Ford iteration~\citep{ford1956network,bellman1958routing}, which solves an instance of the algebraic path problem~\citep{lehmann1977algebraic,baras2010path}. \citet[Chapter~2]{fong2019invitation} give the enriched-category perspective on shortest paths. The bisheaf formulation replaces the scalar stalks $\Quant^{\op}$ of this iteration by heterogeneous complete $\quant$-categories and its identity maps by general adjoint pairs of restrictions and coextensions.

    Fix a vertex $r$ and let $x[0]_r = 0$ and $x[0]_v = \infty$ for $v \neq r$. By \Cref{cor:kleene-powers}, $x[t]_v$ is the length of a shortest walk from $v$ to $r$ with at most $t$ edges, so the flow stabilizes after at most $|V| - 1$ steps. By \Cref{prop:kleene}(b), its limit, the distance to $r$, is the greatest $1$-Lipschitz $0$-cochain below $x[0]$. For explicit values, let $V = \{r,a,b,c\}$, with edges $ra$, $rb$, $ab$, $bc$, and $ac$ of lengths $1$, $4$, $2$, $1$, and $5$. In the order $(r,a,b,c)$, the flow passes through $(0,1,4,\infty)$ and $(0,1,3,5)$ and stops at the distances $(0,1,3,4)$ (\Cref{fig:shortest-paths}).
\end{example}

\begin{figure}[tp]
  \centering
  \begin{tikzpicture}[
      x=1cm, y=1cm,
      vtx/.style={circle, draw=black, fill=white, line width=0.6pt,
                  minimum size=5.2mm, inner sep=0pt, font=\footnotesize},
      rootv/.style={vtx, fill=black!18, draw=black, line width=0.75pt},
      updated/.style={vtx, draw=black, double, double distance=1.0pt, line width=0.55pt},
      gedge/.style={black!45, line width=0.6pt},
      relaxed/.style={black, line width=1.25pt,
                    -{Straight Barb[length=3.5pt,width=4pt]}},
      treeedge/.style={black!22, line width=4.2pt, line cap=round, line join=round},
      lenl/.style={font=\scriptsize, fill=white, inner sep=1pt, text=black!75},
      vname/.style={font=\footnotesize\itshape, text=black},
      panell/.style={font=\small},
      vecl/.style={font=\footnotesize, anchor=north},
    ]

    \newcommand{\sppanel}[9]{%
      \begin{scope}[xshift=#1]
        \coordinate (r) at (0,0.8);
        \coordinate (a) at (1.4,1.6);
        \coordinate (b) at (1.4,0);
        \coordinate (c) at (2.8,0.8);
        \node[panell] at (0.1,2.15) {#2};
        \draw[gedge] (r) -- (a)  (r) -- (b)  (a) -- (b)  (b) -- (c)  (a) -- (c);
        #8
        \path (r) -- (a) node[lenl, midway] {$1$};
        \path (r) -- (b) node[lenl, midway] {$4$};
        \path (a) -- (b) node[lenl, midway] {$2$};
        \path (b) -- (c) node[lenl, midway] {$1$};
        \path (a) -- (c) node[lenl, midway] {$5$};
        \node[rootv] (R) at (r) {$#3$};
        #7
        \node[vname, above left=0pt of R] {$r$};
        \node[vname, above=1pt of A] {$a$};
        \node[vname, below left=0pt of B] {$b$};
        \node[vname, above right=0pt of C] {$c$};
        \node[vecl] at (1.4,-0.7) {#9};
      \end{scope}}

    \sppanel{0cm}{$t=0$}{0}{\infty}{\infty}{\infty}{%
      \node[vtx] (A) at (a) {$\infty$};
      \node[vtx] (B) at (b) {$\infty$};
      \node[vtx] (C) at (c) {$\infty$};}{}%
      {$x[0] = (0,\infty,\infty,\infty)$}

    \sppanel{4.1cm}{$t=1$}{0}{1}{4}{\infty}{%
      \node[updated] (A) at (a) {$1$};
      \node[updated] (B) at (b) {$4$};
      \node[vtx] (C) at (c) {$\infty$};}{%
      \draw[relaxed] (r) -- ($(r)!0.72!(a)$);
      \draw[relaxed] (r) -- ($(r)!0.72!(b)$);}%
      {$x[1] = (0,1,4,\infty)$}

    \sppanel{8.2cm}{$t=2$}{0}{1}{3}{5}{%
      \node[vtx] (A) at (a) {$1$};
      \node[updated] (B) at (b) {$3$};
      \node[updated] (C) at (c) {$5$};}{%
      \draw[relaxed] (a) -- ($(a)!0.72!(b)$);
      \draw[relaxed] (b) -- ($(b)!0.72!(c)$);}%
      {$x[2] = (0,1,3,5)$}

    \sppanel{12.3cm}{$t=3$}{0}{1}{3}{4}{%
      \node[vtx] (A) at (a) {$1$};
      \node[vtx] (B) at (b) {$3$};
      \node[updated] (C) at (c) {$4$};}{%
      \draw[treeedge] (r) -- (a) -- (b) -- (c);
      \draw[gedge] (r) -- (a) -- (b) -- (c);
      \draw[relaxed] (b) -- ($(b)!0.72!(c)$);}%
      {$x[3] = (0,1,3,4)$}

  \end{tikzpicture}
  \caption{Harmonic flow of \Cref{ex:shortest-paths} on the graph with vertices $r$, $a$, $b$, and $c$ and edge lengths $\ell(ra) = 1$, $\ell(rb) = 4$, $\ell(ab) = 2$, $\ell(bc) = 1$, and $\ell(ac) = 5$. Each vertex carries the value $x[t]_v$ of the $0$-cochain at step $t$, the root $r$ (shaded) stays at $0$, and a vertex whose value decreases at step $t$ is marked with a double circle, with the edge $e$ along which the minimum $\ell(e) + x[t-1]_w$ is attained drawn as a bold arrow from the neighbor $w$. At step $t$ the value $x[t]_v$ is the length of a shortest walk from $v$ to $r$ with at most $t$ edges: the estimate $x[1]_b = 4$ along the edge $rb$ is improved at $t = 2$ by the two-edge walk through $a$, and the estimate $x[2]_c = 5$ is improved at $t = 3$ by the three-edge walk through $b$ and $a$. The flow stops at $t = 3 = |V| - 1$, where the last arrows drawn into the vertices $a$, $b$, and $c$ assemble into the shortest-path tree rooted at $r$ (highlighted in gray), and the limit $(0,1,3,4)$ is the greatest $1$-Lipschitz $0$-cochain below $x[0]$.}
  \label{fig:shortest-paths}
\end{figure}

A second example glues local estimates over a covering. It is a discrete instance of the bisheaves of presheaves of \Cref{ex:powers}, over $\law$, and its harmonic flow is a path problem on the nerve of the covering.

\begin{figure}[tp]
    \centering
    \begin{tikzpicture}[scale=1.15, every node/.style={font=\small}]
        \begin{scope}
            \node at (-0.9,2.75) {(a)};
            \fill[black!7, rounded corners=10pt]  (-0.45,0.55)  rectangle (1.45,2.45);
            \fill[black!11, fill opacity=0.55, rounded corners=10pt] (0.55,0.55)  rectangle (2.45,2.45);
            \fill[black!11, fill opacity=0.55, rounded corners=10pt] (-0.45,-0.45) rectangle (1.45,1.45);
            \fill[black!7, fill opacity=0.55, rounded corners=10pt]  (0.55,-0.45) rectangle (2.45,1.45);
            \draw[black!60, thick, rounded corners=10pt] (-0.45,0.55)  rectangle (1.45,2.45);
            \draw[black!65, thick, densely dashed, rounded corners=10pt] (0.55,0.55)  rectangle (2.45,2.45);
            \draw[black!65, thick, densely dashed, rounded corners=10pt] (-0.45,-0.45) rectangle (1.45,1.45);
            \draw[black!60, thick, rounded corners=10pt] (0.55,-0.45) rectangle (2.45,1.45);
            \node[black!85] at (-0.2,2.25)  {$U_1$};
            \node[black!85] at (2.2,2.25)   {$U_2$};
            \node[black!85] at (-0.2,-0.25) {$U_3$};
            \node[black!85] at (2.2,-0.25)  {$U_4$};
            \foreach \x in {0,1,2} \foreach \y in {0,1} \draw[black!45] (\x,\y) -- (\x,\y+1);
            \foreach \x in {0,1} \foreach \y in {0,1,2} \draw[black!45] (\x,\y) -- (\x+1,\y);
            \draw[black!45] (0,2) -- (1,1) -- (2,0);
            \draw[black!45] (1,1) -- (2,2);
            \foreach \x in {0,1,2} \foreach \y in {0,1,2} \fill[black] (\x,\y) circle (1.7pt);
            \coordinate (n1) at (0.5,1.5);
            \coordinate (n2) at (1.5,1.5);
            \coordinate (n3) at (0.5,0.5);
            \coordinate (n4) at (1.5,0.5);
            \draw[black, very thick] (n1) -- (n2) -- (n4) -- (n3) -- (n1);
            \draw[black, very thick] (n1) to[bend left=25] (n4);
            \draw[black, very thick] (n3) to[bend left=25] (n2);
            \foreach \n in {n1,n2,n3,n4} \filldraw[fill=white, draw=black, line width=1.1pt] (\n) circle (2.7pt);
        \end{scope}
        \begin{scope}[xshift=4.6cm, yshift=-0.1cm]
            \node at (-0.7,2.85) {(b)};
            \fill[black!7] plot[smooth cycle, tension=0.7] coordinates {(-0.5,2.5) (1.3,2.4) (3.0,1.7) (2.9,1.0) (2.0,0.15) (1.0,-0.3) (0.1,-0.3) (-0.5,0.8)};
            \fill[black!12, fill opacity=0.55] plot[smooth cycle, tension=0.7] coordinates {(0.1,0.6) (0.3,-0.4) (1.5,-0.75) (2.7,-0.75) (3.3,-0.1) (3.3,1.3) (2.8,2.0) (2.1,1.75) (1.25,1.2)};
            \draw[black!65, thick] plot[smooth cycle, tension=0.7] coordinates {(-0.5,2.5) (1.3,2.4) (3.0,1.7) (2.9,1.0) (2.0,0.15) (1.0,-0.3) (0.1,-0.3) (-0.5,0.8)};
            \draw[black!65, thick, densely dashed] plot[smooth cycle, tension=0.7] coordinates {(0.1,0.6) (0.3,-0.4) (1.5,-0.75) (2.7,-0.75) (3.3,-0.1) (3.3,1.3) (2.8,2.0) (2.1,1.75) (1.25,1.2)};
            \node[black] at (-0.85,1.4) {$U_i$};
            \node[black] at (3.5,-0.6)  {$U_j$};
            \coordinate (a) at (0,2);
            \coordinate (b) at (1,1.5);
            \coordinate (c) at (2.6,1.4);
            \coordinate (d) at (0.6,0.2);
            \coordinate (f) at (2.4,-0.2);
            \draw[black!50, -latex] (a) -- (b);
            \draw[black!50, -latex] (b) -- (c);
            \draw[black!50, -latex] (a) -- (d);
            \draw[black!50, -latex] (c) -- (f);
            \draw[black!50, -latex] (d) -- (f);
            \draw[black, very thick, -latex] (d) -- (c)
                node[midway, above=1pt, sloped, black] {$A_i(d,c) = 5$}
                node[midway, below=1pt, sloped, black] {$A_j(d,c) = 8$};
            \foreach \p in {a,b,c,d,f} \fill[black] (\p) circle (1.8pt);
            \node[above left] at (a) {$a$};
            \node[above] at (b) {$b$};
            \node[above right] at (c) {$c$};
            \node[below left] at (d) {$d$};
            \node[below right] at (f) {$f$};
        \end{scope}
    \end{tikzpicture}
    \caption{Reconciliation of local cost estimates over the nerve of a covering (\Cref{ex:map-merging}). (a)~Network of waypoints (solid vertices, thin edges) covered by four regions, and the $1$-skeleton of the nerve of the covering (open vertices, thick edges), which is the graph $G$ of \Cref{ex:map-merging}. (b)~Two regions $U_i = \{a,b,c,d\}$ (solid boundary) and $U_j = \{c,d,f\}$ (dashed boundary) with overlap $U_{ij} = \{c,d\}$, whose agents disagree about the cost of the shared route $d \to c$. At the tolerance $2$ in both directions, one step of the harmonic flow raises the estimate of region $i$ to $6$ and leaves that of region $j$ at $8$, after which the two estimates agree up to the tolerance. Had the overlap consisted of the single waypoint $c$, the only shared quantity would be the cost $A(c,c)$ of the loop at $c$, and the regions would share no route between distinct waypoints.}
    \label{fig:waypoints}
\end{figure}

\begin{example}[Reconciling local cost estimates]\label{ex:map-merging}
    Let $X$ be a finite set of \emph{waypoints} covered by \emph{regions} $U_1, \dots, U_n \subseteq X$.
    An agent surveying $U_i$ records its estimate of the cost of traveling directly between waypoints as a matrix $A_i : U_i \times U_i \to [0,\infty]$, with $A_i(a,b) = \infty$ when it judges that there is no direct route from $a$ to $b$.
    Writing $z_{+} \coloneqq \max(z,0)$, with $(\infty - \infty)_{+} \coloneqq 0$ as in \Cref{ex:quantales}(b), the estimates on a finite set $U$ form the power $\Quant^{U \times U}$ over $\quant = \law$, complete by \Cref{lem:self-enrichment}, with
    \[ \hom(A,B) = \max_{a,b \in U} \bigl( B(a,b) - A(a,b) \bigr)_{+}, \]
    so that $A \qleq_q B$ if and only if $B \leq A + q$ pointwise.
    The underlying order is the reversed pointwise order, conical meets are pointwise maxima, cotensors are $q \cotensor A = (A - q)_{+}$, and the top object is the zero matrix.

    Let $G$ be the $1$-skeleton of the nerve of the covering, with the regions as vertices and an edge $\{i,j\}$ whenever the overlap $U_{ij} \coloneqq U_i \cap U_j$ is nonempty (\Cref{fig:waypoints}).
    The sets $Y_i \coloneqq U_i \times U_i$ and $Y_e \coloneqq U_{ij} \times U_{ij}$ for $\partial e = \{i,j\}$, with the inclusions $Y_e \hookrightarrow Y_i$, determine a biadjoint bisheaf of powers (\Cref{ex:powers}): the restriction $\ushf_{i \face e}$ restricts a matrix to the rows and columns indexed by $U_{ij}$, and the coextension $\oshf_{i \face e}$ extends a matrix on $U_{ij}$ by $0$, the top element of $\law$, outside the overlap.
    Read the tolerance $W(i,j) \in [0,\infty]$ as the number of cost units by which region $i$ permits region $j$ to overestimate a shared cost.
    By \Cref{thm:global-sections}, a family $A = (A_i)_{i \in V}$ of local estimates is a weighted global section if and only if
    \begin{equation*}
        A_j(a,b) \leq A_i(a,b) + W(i,j) \qquad \text{for every edge } e, \text{ both orderings } (i,j) \text{ of } \partial e, \text{ and all } a,b \in U_{ij},
    \end{equation*}
    that is, if neighboring regions agree on the cost of every shared route up to the tolerance.
    Since cotensors subtract and conical meets are maxima, the harmonic flow is the distributed update
    \[ A_i[t+1](a,b) = \max\Bigl( A_i[t](a,b), \max_{\substack{j \sim_e i \\ a,b \in U_{ij}}} \bigl( A_j[t](a,b) - W(i,j) \bigr)_{+} \Bigr), \]
    in which every region raises its estimate of each shared route to within tolerance of the most pessimistic estimate held by a neighbor.
    For a walk $\gamma = (i_0, i_1, \dots, i_s)$ in $G$, write $\ell_W(\gamma) \coloneqq \sum_{r < s} W(i_{r+1}, i_r)$ for the tolerance accumulated along $\gamma$. The bisheaf is biadjoint, so \Cref{cor:kleene-powers} and induction on $t$ give
    \[ A_i[t](a,b) = \max_{\gamma} \bigl( A_{j}[0](a,b) - \ell_W(\gamma) \bigr)_{+}, \]
    where the maximum ranges over the walks $\gamma$ of length $s \leq t$ from a region $j = i_0$ to $i_s = i$ that pass only through regions containing both $a$ and $b$, the walk of length $0$ included.
    The maximum is attained on simple paths, so the flow stabilizes after at most $n - 1$ steps, and by \Cref{prop:kleene}(b) its limit is the greatest weighted global section below $A[0]$: the cheapest family of mutually consistent estimates that undercuts no initial estimate.
    The flow is thus a path problem on the nerve in the spirit of \Cref{ex:shortest-paths}, and \Cref{fig:waypoints}(b) gives explicit values.
\end{example}

The \emph{synchronization problem} asks for values at the vertices of a graph that realize relations prescribed along its edges, as in the synchronization $x_w = g_e \scdot x_v$ of orthogonal transformations $g_e \in O(d)$ that underlies the graph connection Laplacian of \citet{singer2012}, and in its variants over permutation groups~\citep{pachauri2013} and principal bundles~\citep{gao_geometry_2021}.
In the quantale analogue, elements $\tau_e \in \quant$ are prescribed along oriented edges $e : v \to w$, the problem asks for $0$-cochains with $x_w = x_v \odot \tau_e$, exactly or up to a tolerance, and residuation takes the place of inversion.

\begin{example}[Quantale synchronization]\label{ex:sync}
    Synchronization requires an orientation of each edge, which the undirected graphs of \Cref{sec:sheaves} do not supply, and we fix an orientation as additional data, leaving the incidence poset unchanged. An \emph{orientation} of $G$ chooses for every edge $e$ a source and a target among its endpoints, and we write $e : v \to w$ for an edge with source $v$ and target $w$. A \emph{$\quant$-connection} on an oriented graph $G$ is a map $\tau : E \to \quant$ with values $\tau_e$, and the \emph{synchronization problem} for $(G, \tau)$ asks for the tuples $x \in \prod_{v \in V} \quant$ with
    \[ x_w = x_v \odot \tau_e \qquad \text{for every edge } e : v \to w. \]

    Let $\tau$ be a $\quant$-connection on an oriented graph $G$ and $W$ a weight function on $G$. Let every stalk be $\Quant$, and for every edge $e : v \to w$ let
    \[ \ushf_{v \face e} \coloneqq (-) \odot \tau_e, \qquad \ushf_{w \face e} \coloneqq \id, \qquad \oshf_{v \face e} \coloneqq [\tau_e, -], \qquad \oshf_{w \face e} \coloneqq \id. \]
    These maps form an adjoint bisheaf, the \emph{synchronization bisheaf}. Indeed, the inequality $[p,p'] \odot p \odot \tau_e \preceq p' \odot \tau_e$ gives $[p,p'] \preceq [p \odot \tau_e, p' \odot \tau_e]$, and \Cref{lem:functorality-of-hom}(7) gives $[p,p'] \preceq \bigl[[\tau_e,p], [\tau_e,p']\bigr]$, so restriction and coextension are $\quant$-functors. Moreover, the adjunction at $v$ is the identity $[x \odot \tau_e, z] = [x, [\tau_e, z]]$ of \Cref{lem:functorality-of-hom}(6).

    The stalk $\Quant$ is complete by \Cref{lem:self-enrichment}(a), with cotensors $\dlhom{q}{-}$ and conical meets the meets of $\quant$, so the Lawvere Laplacian $L^W$ exists. The transports into $v$ are $x \mapsto x \odot \tau_e$ along edges $e : u \to v$ and $x \mapsto [\tau_e, x]$ along edges $e : v \to w$, and the decomposition~\eqref{eq:weighted-meet-formula} gives
    \begin{equation}\label{eq:sync-lap}
        (L^W x)_v \qequals_1 \bigmeet_{e : u \to v} \dlhom{W(v, u)}{x_u \odot \tau_e} \meet \bigmeet_{e : v \to w} \dlhom{W(v, w)}{[\tau_e, x_w]}.
    \end{equation}
    For an edge $e : v \to w$, the condition of \Cref{thm:global-sections} for the ordered pair $(v,w)$ is $W(v,w) \preceq [x_v \odot \tau_e, x_w]$, and for $(w,v)$ it is $W(w,v) \preceq [x_w, x_v \odot \tau_e]$. By the residuation~\eqref{eq:adjoint-hom}, a $0$-cochain $x$ is a weighted global section, equivalently a suffix point of $L^W$ by \Cref{thm:hodge-thm}, if and only if, for every edge $e : v \to w$,
    \begin{equation}\label{eq:sync-sections}
        W(v, w) \odot x_v \odot \tau_e \preceq x_w \qquad \text{and} \qquad W(w, v) \odot x_w \preceq x_v \odot \tau_e.
    \end{equation}
    At $W \equiv 1$ the two inequalities combine to $x_w = x_v \odot \tau_e$, and the weighted global sections are the solutions of the synchronization problem.

    Let $C$ be a directed cycle $v_0 \to v_1 \to \dots \to v_n = v_0$ with edges $e_1, \dots, e_n$, and let its \emph{holonomy} be the product $h_C \coloneqq \tau_{e_1} \odot \dots \odot \tau_{e_n}$. Along $C$ the synchronization equations compose to $x_{v_0} = x_{v_0} \odot h_C$. Over a group, a solution exists only if the holonomy of every cycle is trivial. Quantale multiplication is not cancellative, and the bottom $0$-cochain $x \equiv \bot$ is a weighted global section of every synchronization bisheaf for every weight function, the analogue of the zero section of a linear sheaf. Heuristically, inconsistent cycles drive solutions toward absorbing values instead of rendering the problem infeasible. For instance, let $M$ be a positive integer and consider the finite quantale $\{0,1,\dots,M\}$, ordered by $\geq$, with multiplication $(a,b) \mapsto \min(a+b,M)$ and unit $0$, a truncation of $\law$. On a directed $3$-cycle $u \to v \to w \to u$ with $\tau_e = 1$ on every edge (the integer $1$, not the unit $0$) and the unit weight function $W \equiv 0$, the harmonic flow from $x[0]=(0,0,0)$ increments all three entries by $1$ at each step and reaches the absorbing bottom $0$-cochain $(M,M,M)$ after $M$ steps. Under the hypothesis of \Cref{prop:kleene}(a), which holds for instance when $\quant$ is finite (\Cref{rem:chain-preservation}), the harmonic flow from initial estimates $x[0]$ computes the greatest weighted global section below them.
\end{example}

\begin{example}[Clocks]\label{ex:clocks}
    Over the extended Lawvere quantale $\elaw$ of \Cref{ex:quantales}(e), the synchronization equations read $x_w = x_v + \tau_e$ for $e : v \to w$, which is clock synchronization from pairwise offset measurements, with $\tau_e \in \mathbb{R}$ the measured offset of $w$ relative to $v$. For finite values, the weighted global sections are the solutions of the synchronization bounds~\eqref{eq:sync-sections}, which read
    \[ x_w - (x_v + \tau_e) \leq W(v, w) \qquad \text{and} \qquad (x_v + \tau_e) - x_w \leq W(w, v), \]
    the two tolerances bounding the error in each direction separately. The tolerance $W(w, v) = \bot = \infty$ leaves only the first constraint, so that $x_w$ may fall arbitrarily far below $x_v + \tau_e$ but may exceed it by at most $W(v, w)$. A negative tolerance $W(v,w) = -\delta < 0$ lies strictly above the unit $0$ of $\elaw$, and \Cref{def:weight-function} then requires $W(w,v) \geq \delta$, in accordance with \Cref{lem:above-unit}. The first constraint becomes the margin $(x_v + \tau_e) - x_w \geq \delta$, so that $x_w$ falls below $x_v + \tau_e$ by at least $\delta$. Over $\law$ the offsets would have to be nonnegative, and the coextension $[\tau_e, -]$ would be truncated subtraction. We therefore work over $\elaw$.

    For explicit values, take the path $u \to v \to w$ with offsets $2$ and $3$, the tolerance $1$ on every ordered pair of adjacent vertices (one unit of time, not the unit $0$ of $\elaw$), and initial estimates $x[0] = (0,0,0)$. The order of $\elaw$ reverses the numerical order, so meets are numerical maxima, and by the Laplacian~\eqref{eq:sync-lap} the harmonic flow replaces each $x_v$ by the largest of $x_v$, the values $x_u + \tau_e - 1$ for edges $e : u \to v$, and the values $x_w - \tau_e - 1$ for edges $e : v \to w$. It passes through $(0,1,2)$ and stops at $(0,1,3)$, where both edges are at the limit of the tolerance: $(x_u + 2) - x_v = (x_v + 3) - x_w = 1$. The flow stabilizes, so by \Cref{prop:kleene}(b) the limit $(0,1,3)$ is the greatest weighted global section below $x[0]$ in the order of $\elaw$, that is, the componentwise least solution above $(0,0,0)$ in the numerical order.
\end{example}

\begin{example}[Synchronization over binary relations]\label{ex:matchings}
    Let $\quant = \cat{Rel}(X)$ be the quantale of binary relations on $X = \{0,1\}$ (\Cref{ex:quantales}(c)), which is neither commutative nor affine. Its right hom is $[S,T] = \{(i,j) : (i,k) \in T \text{ for all } k \in X \text{ with } (j,k) \in S\}$, and its left hom is $\dlhom{S}{T} = \{(j,k) : (i,k) \in T \text{ for all } i \in X \text{ with } (i,j) \in S\}$. Consider the synchronization bisheaf of \Cref{ex:sync} on a single oriented edge $e : u \to v$ with $\tau_e = \{(0,1)\}$, and let $W(u,v) = \Delta_X \cup \{(1,0)\}$ and $W(v,u) = \emptyset$. These tolerances form a weight function, since both products $W(u,v) \odot W(v,u)$ and $W(v,u) \odot W(u,v)$ are empty. The tolerance $W(u,v)$ lies strictly above the unit $\Delta_X$ and is not central.

    For a $0$-cochain $x = (x_u, x_v)$, the Laplacian~\eqref{eq:sync-lap} reads $(L^W x)_u = \dlhom{W(u,v)}{[\tau_e, x_v]}$ and $(L^W x)_v = X \times X$, the top element. The map $L^W$ is monotone (\Cref{lem:lap-functor}), but it is not a $\quant$-functor. For $x = (X \times X, \{(1,1)\})$ and $y = (X \times X, \{(0,1)\})$, the pair $(0,1)$ lies in $\hom(x,y)$ and not in $\hom(L^W x, L^W y)$. The Hodge--Lawvere theorem (\Cref{thm:hodge-thm}) and the convergence of the harmonic flow (\Cref{prop:kleene}) do not require central weight functions, and they apply here, whereas \Cref{prop:gamma-limits} does not. Its conclusions fail as well. For $q = \{(0,0)\}$ and the weighted global section $x = (\{(0,0)\}, \{(0,1),(1,1)\})$, the componentwise tensor $(\{(0,0)\}, \{(0,1)\})$ is not a weighted global section, and the tensor of $x$ by $q$ does not exist in $\wsections$, by exhaustive enumeration. Consequently, $\wsections$ is neither complete nor coreflective, although it is closed under conical joins and contains a greatest weighted global section below every $0$-cochain. From $x[0] = (\{(1,0)\}, \{(1,0)\})$, we compute $[\tau_e, x[0]_v] = \{(0,1),(1,1)\}$ and $(L^W x[0])_u = \{(0,1),(1,1)\}$. The harmonic flow therefore reaches $x[1] = (\emptyset, \{(1,0)\})$ in one step and stops there, at the greatest weighted global section below $x[0]$.
\end{example}

\medskip\noindent Discrete event systems supply a bisheaf of a different kind, whose stalks are products of $\Quant^{\op}$ indexed by a set of events and whose restrictions are matrices of delays.

\begin{example}[Discrete event systems]\label{ex:maxplus}
    Discrete event systems, whose state advances at the occurrence of discrete events, resist linear modeling over a field, but a large class of such systems is linear over the max-plus semiring~\citep{butkovic2010,deshutter2020}.
    Fix a finite set $I$ of \emph{events}, and let $\quant = \elaw$ be the extended Lawvere quantale of \Cref{ex:quantales}(e). Because $\elaw$ is ordered by $\geq$, its joins are numerical infima and it is a min-plus quantale, whereas passing to the opposite $\elaw$-category $(\Quant^{I})^{\op}$ of \emph{timing vectors} $x : I \to [-\infty,\infty]$, where $x(i)$ is the epoch at which event $i$ occurs, reverses the underlying order to the componentwise numerical order $\leq$.
    This stalk is complete by \Cref{lem:self-enrichment}(b) and~(c): its conical meets are componentwise minima, its conical joins are componentwise maxima, cotensoring by $q$ adds $q$, and tensoring by $q$ subtracts $q$.
    A max-plus system is specified by a matrix $A : I \times I \to [0,\infty)$ of \emph{delays} and, for each $j \in I$, a nonempty set $\pre(j) \subseteq I$ of events on which $j$ depends, and it advances by the map
    \begin{equation}\label{eq:max-plus}
        T(x) = \Bigl( \textstyle\max_{i \in \pre(j)} \bigl( x(i) + A(i,j) \bigr) \Bigr)_{j \in I}.
    \end{equation}
    The map $T$ preserves componentwise maxima and commutes with tensors, so it is a $\quant$-functor that preserves weighted joins by the decomposition~\eqref{eq:weighted-meet-formula}, and \Cref{prop:kan-bisheaf} gives its right adjoint $T_{\ast}$, written $T^{\sharp}$ in max-plus algebra~\citep{butkovic2010},
    \[ T^{\sharp}(y)(i) = \min_{j : i \in \pre(j)} \bigl[A(i,j), y(j)\bigr], \]
    which equals $\min_{j : i \in \pre(j)} \bigl(y(j) - A(i,j)\bigr)$ for finite values, the minimum over the empty family being $\infty$.
    Over $\law$, tensoring is truncated subtraction, which $T$ does not preserve, so we work over $\elaw$.

    Coupled systems are glued along a finite simple graph $G$ whose vertices index subsystems and whose edges index pairwise couplings. The stalk at a cell $\sigma$ is $(\Quant^{I_\sigma})^{\op}$ for a finite set $I_\sigma$ of events, the restrictions are max-plus maps $T_{v \face e} : (\Quant^{I_v})^{\op} \to (\Quant^{I_e})^{\op}$ of the form~\eqref{eq:max-plus} with delay matrices $A_{v \face e} : I_v \times I_e \to [0,\infty)$ and nonempty predecessor sets $\pre_{v \face e}(j) \subseteq I_v$, and the coextensions are their residuals.
    These data form an adjoint bisheaf, and its Lawvere Laplacian is
    \[ (L^W x)_v(i) = \min_{w \sim_e v} \Bigl( W(v,w) + T_{v \face e}^{\sharp}\bigl(T_{w \face e} x_w\bigr)(i) \Bigr). \]
    The tolerances $W(v,w) = W(w,v) = 0$ demand exact synchrony across the edge, and $W(v,w) = \infty$ imposes no constraint in the direction from $v$ to $w$.
    By \Cref{thm:hodge-thm}, the timing assignments synchronized up to $W$ are the suffix points of $L^W$. Since the sets $I_\sigma$ and the graph $G$ are finite, $L^W$ is a composite of finitely many maxima, minima, and translations by elements of $[-\infty,\infty]$, each of which preserves the infima of descending chains in $[-\infty,\infty]$. Consequently, $L^W$ preserves such infima, and by \Cref{prop:kleene} the harmonic flow computes the greatest synchronized assignment below the initial timings.
    For symmetric weight functions, $L^W$ has the form of the tropical Tarski Laplacian of \citet[Definition~2]{riess2023}, which is stated for square matrices with a finite entry in every row and column and for finite tolerances on the edges, and the harmonic flow has the form of their heat equation~\citep[Eq.~(9)]{riess2023}. By \Cref{thm:hodge-thm}, a $0$-cochain with finite values is a fixed point of that heat equation if and only if $\norm{T_{v \face e} x_v - T_{w \face e} x_w}_{\infty} \leq W(v,w)$ for every edge $e$ with $\partial e = \{v,w\}$. The forward implication is \citet[Theorem~3]{riess2023}.

    For explicit values, consider three workcells $u$, $v$, and $w$ arranged in the chain $u \to v \to w$, each with event set $I_c = \{s,f\}$ consisting of a start event $s$ and a finish event $f$ and with processing times $p_u = 3$, $p_v = 2$, and $p_w = 4$. Consecutive cells $c \to c'$ are coupled through a buffer $e$ with a single part-arrival event $I_e = \{a\}$. The upstream restriction $T_{c \face e}(x_c) = \max\bigl(x_c(s) + p_c + 1, x_c(f) + 1\bigr)$ has $\pre_{c \face e}(a) = \{s,f\}$ with delays $A_{c \face e}(s,a) = p_c + 1$ and $A_{c \face e}(f,a) = 1$, which include one time unit of transfer to the buffer, and the downstream restriction $T_{c' \face e}(x_{c'}) = x_{c'}(s)$ has $\pre_{c' \face e}(a) = \{s\}$ with delay $A_{c' \face e}(s,a) = 0$. The upstream restriction places the arrival one time unit after the later of the recorded finish $x_c(f)$ and the earliest possible finish $x_c(s) + p_c$. We set the tolerances to $W(u,v) = W(v,w) = 0$ toward each downstream cell and $W(v,u) = W(w,v) = 2$ toward each upstream cell, so a schedule is a weighted global section if and only if each downstream cell starts no earlier than its part arrives and at most two time units later. Starting from the initial schedule $x[0]_u = (1,4)$, $x[0]_v = (9,11)$, and $x[0]_w = (7,11)$, the harmonic flow lowers $v$ to $x[1]_v = (4,6)$ at the first step and $u$ to $x[2]_u = (0,3)$ at the second, stabilizing at $x[2] = \bigl((0,3), (4,6), (7,11)\bigr)$, the latest synchronized schedule below the initial schedule. The flow preserves the processing constraint $x_c(f) \geq x_c(s) + p_c$, which the initial schedule satisfies.
\end{example}


\section{Conclusion}\label{sec:conclusion}

Cellular sheaves of $\quant$-categories and the Lawvere Laplacian provide a categorification of harmonic analysis on graphs.
When the vertex stalks are complete and the restrictions have right adjoints, the weighted global sections $[\widetilde{W},\shf]$ coincide with the suffix points of the Lawvere Laplacian (\Cref{thm:hodge-thm}), to which the harmonic flow converges under a continuity hypothesis (\Cref{prop:kleene}).
For central weight functions, the flow computes the coreflection of the $0$-cochains onto the weighted global sections. Several directions remain open:
\begin{itemize}
    \item[(a)] Cellular $\quant$-sheaves are defined over every finite preorder, such as the simplex poset of a semisimplicial set~\citep[Observations~2.7 and~4.16]{fioresi2026}, whereas higher Laplacians require a quantale-valued analogue of the cochain complex. Likewise, sheaves on directed graphs, indexed by a category of elements~\citep{althausCompositionalAlgorithmsCompositional2023}, would make the edge orientations of \Cref{ex:sync} intrinsic.
    \item[(b)] Allowing the base quantale to vary across cells via the Grothendieck construction~\citep{moeller2020monoidal} points toward a Hodge--Lawvere theorem over a quantaloid~\citep{stubbe_introduction_2014}, which waits in the wings. Beyond thin bases, enriching in a closed monoidal category~\citep{Kelly2005} or bicategory would replace posetal inequalities and meets with morphisms and ends.
    \item[(c)] Because quantale joins are idempotent and $\bot$ absorbs under $\odot$, the Lawvere Laplacian does not subsume linear Laplacians. Whether the Chu construction~\citep{barr1996}, which lifts linear and enriched adjoints simultaneously and underlies conic duality over complete $\quant$-categories~\citep[Part~IV]{lopezEnrichedCategoriesQuantales2025}, connects the two operators remains open.
\end{itemize}

\appendix


\section{Residuation and Weighted Meets and Joins}\label{appx:technical}

We use the following standard properties of the two homs of a unital quantale. Over a commutative base the two homs coincide, and the right-hom and left-hom forms of each clause agree.

\begin{lemma}\label{lem:functorality-of-hom}
    Let $\quant$ be a unital quantale, let $p,q,r \in \quant$, and let $(q_i)_{i \in I}$ be a family in $\quant$. Then the following properties hold, and clauses~(1)--(3) also hold with $\dlhom{-}{-}$ in place of $[-,-]$.
    \begin{enumerate}[label=(\arabic*)]
        \item If $q \preceq r$, then $[p, q] \preceq [p, r]$. If $p \preceq r$, then $[r, q] \preceq [p, q]$.
        \item We have $\bigl[p, \bigmeet_{i \in I} q_i\bigr] = \bigmeet_{i \in I} [p, q_i]$ and $\bigl[\bigjoin_{i \in I} q_i, p\bigr] = \bigmeet_{i \in I} [q_i, p]$.
        \item The unit satisfies $[1, q] = q$.
        \item The inequality $q \preceq p$ holds if and only if $1 \preceq [q,p]$, if and only if $1 \preceq \dlhom{q}{p}$. In particular, if $\quant$ is affine and $q \preceq p$, then $[q, p] = \dlhom{q}{p} = 1$.
        \item We have $r \odot [p, q] \preceq [p, r \odot q]$ and $\dlhom{p}{q} \odot r \preceq \dlhom{p}{q \odot r}$.
        \item We have $[p, [q ,r]] = [p \odot q , r]$ and $\dlhom{p}{\dlhom{q}{r}} = \dlhom{q \odot p}{r}$, as well as $\dlhom{p}{[q,r]} = [q, \dlhom{p}{r}]$. Over a commutative base, $[p,[q,r]] = [q,[p,r]]$, and this identity for all $p,q,r$ is equivalent to commutativity.
        \item We have $[q,r] \odot [p,q] \preceq [p,r]$ and $\dlhom{p}{q} \odot \dlhom{q}{r} \preceq \dlhom{p}{r}$.
    \end{enumerate}
\end{lemma}

\begin{proof}
    Each clause follows from the residuation~\eqref{eq:adjoint-hom} by comparing lower sets. For instance, in clause~(6), $c \preceq \dlhom{p}{[q,r]}$ if and only if $p \odot c \odot q \preceq r$, if and only if $c \preceq [q, \dlhom{p}{r}]$, and in~(7), $[q,r] \odot [p,q] \odot p \preceq [q,r] \odot q \preceq r$. For the last claim of~(6), taking $r = q \odot p$ in $[p \odot q, r] = [q \odot p, r]$ gives $1 \preceq [p \odot q, q \odot p]$, so $p \odot q \preceq q \odot p$ for all $p$ and $q$.
\end{proof}

\medskip\noindent The following lemmas on weighted meets and joins are used in \Cref{sec:prelim,sec:lawvere-laplacian}.

\begin{lemma}\label{lem:weights}
    Let $\C$ be a $\quant$-category, $\S$ a set, $D: \S \to \obj(\C)$ a family of objects, $W: \S \to \quant$ a weight, and $s \in \S$. If $\bigmeet^W D$ exists in $\C$, then $\bigmeet^W D \qleq_{W(s)} D(s)$. If $\bigjoin^W D$ exists in $\C$, then $D(s) \qleq_{W(s)} \bigjoin^W D$.
\end{lemma}

\begin{proof}
    Let $m = \bigmeet^W D$. \Cref{def:weighted} gives
    \begin{align*}
    1 \preceq \hom(m,m) &= \bigmeet_{s'} \dlhom{W(s')}{\hom(m, D(s'))} \\
    &\preceq \dlhom{W(s)}{\hom(m, D(s))},
    \end{align*}
    and transposing gives $W(s) = W(s) \odot 1 \preceq \hom(m, D(s))$. The join half is the same computation with the right hom.
\end{proof}

\begin{lemma}\label{lem:product-hom}
    Let $\C$ be a $\quant$-category, $\S$ a set, $D, D' : \S \to \obj(\C)$ families of objects, and $W, W' : \S \to \quant$ weights with $W'(s) \preceq W(s)$ for every $s \in \S$, such that the weighted meets and joins below exist in $\C$.
    \begin{itemize}
        \item[(a)] If $D(s) \qleq_1 D'(s)$ for every $s$, then $\bigmeet^W D \qleq_1 \bigmeet^{W'} D'$.
        \item[(b)] If the values of $W$ are central, then $\bigmeet_{s \in \S} \hom(D(s), D'(s)) \preceq \hom\bigl(\bigmeet^W D, \bigmeet^{W'} D'\bigr)$.
    \end{itemize}
    Dually, if $D(s) \qleq_1 D'(s)$ for every $s$, then $\bigjoin^{W'} D \qleq_1 \bigjoin^{W} D'$.
\end{lemma}

\begin{proof}
    Let $m = \bigmeet^W D$ and $h_s = \hom(D(s), D'(s))$. By the composition law and \Cref{lem:weights}, we have $h_s \odot W(s) \preceq h_s \odot \hom(m, D(s)) \preceq \hom(m, D'(s))$. (a) If $1 \preceq h_s$, then $W'(s) \preceq W(s) \preceq \hom(m, D'(s))$, so $1 \preceq \dlhom{W'(s)}{\hom(m, D'(s))}$ for every $s$, and $1 \preceq \hom(m, \bigmeet^{W'} D')$ by \Cref{def:weighted}. (b) If $W(s)$ is central, then $W'(s) \odot h_s \preceq W(s) \odot h_s = h_s \odot W(s) \preceq \hom(m, D'(s))$, so $h_s \preceq \dlhom{W'(s)}{\hom(m, D'(s))}$, and taking the meet over $s$ gives the claim by \Cref{def:weighted}. The dual statement follows in the same way from $\hom(D(s), \bigjoin^W D') \succeq \hom(D'(s), \bigjoin^W D') \odot h_s \succeq W(s) \succeq W'(s)$.
\end{proof}

\begin{proof}[Proof of \Cref{lem:self-enrichment}]
    (a) By \Cref{lem:functorality-of-hom}(6), we have $[x, \dlhom{q}{y}] = \dlhom{q}{[x,y]}$ and $[q \odot x, y] = [q, [x,y]]$, which are the defining identities of the cotensor and the tensor in $\Quant$. By \Cref{lem:functorality-of-hom}(2), we have $[x, \bigmeet_i y_i] = \bigmeet_i [x, y_i]$ and $[\bigjoin_i y_i, x] = \bigmeet_i [y_i, x]$, so conical meets and joins are meets and joins in $\quant$. Completeness follows from the decomposition~\eqref{eq:weighted-meet-formula}.

    (b) If $\quant$ is commutative, then $\hom_{\Quant^\op}(x, q \odot y) = [q \odot y, x] = [q, [y,x]]$ and $\hom_{\Quant^\op}([q,y], x) = [x, [q,y]] = [q, [x,y]]$ by \Cref{lem:functorality-of-hom}(6), which are the defining identities of the cotensor and the tensor in $\Quant^\op$. Conical meets and joins in $\Quant^\op$ are joins and meets in $\quant$ by \Cref{lem:functorality-of-hom}(2), and completeness follows as in~(a).

    (c) In a product, we have $\hom(x, (\bigmeet^W D_i)_i) = \bigmeet_i \bigmeet_s \dlhom{W(s)}{\hom(x_i, D(s)_i)} = \bigmeet_s \dlhom{W(s)}{\hom(x, D(s))}$ by \Cref{lem:functorality-of-hom}(2), and dually for weighted joins.
\end{proof}

\begin{proof}[Proof of \Cref{lem:complete-cocomplete}]
    Suppose that $\C$ admits all weighted joins, and let $D : \S \to \obj(\C)$ and $W : \S \to \quant$. Let $\psi(x) \coloneqq \bigmeet_{s \in \S} \dlhom{W(s)}{\hom(x, D(s))}$ for $x \in \C$, and let $m$ be the join of the family $(x)_{x \in \C}$ weighted by $\psi$. We show that $\hom(z, m) = \psi(z)$ for every $z \in \C$, so that $m$ is the weighted meet of $D$ by $W$. By \Cref{lem:weights}, $\psi(z) \preceq \hom(z, m)$. Conversely, for all $s \in \S$ and $x \in \C$, we have $W(s) \odot \psi(x) \preceq \hom(x, D(s))$, that is, $W(s) \preceq [\psi(x), \hom(x, D(s))]$. Taking the meet over $x$ gives $W(s) \preceq \hom(m, D(s))$ by \Cref{def:weighted}, and the composition law gives $W(s) \odot \hom(z, m) \preceq \hom(m, D(s)) \odot \hom(z,m) \preceq \hom(z, D(s))$. Hence $\hom(z,m) \preceq \dlhom{W(s)}{\hom(z, D(s))}$ for every $s$, and $\hom(z, m) \preceq \psi(z)$. The converse is dual.
\end{proof}


\section{Weighted Global Sections as a Weighted Limit}\label{appx:weighted-limit}

The category $\QCat$ is a locally preordered $2$-category whose hom-sets carry the structure of the $\quant$-categories $[\C,\D]$ of \Cref{ex:constructions}(c), and composition $(g,f) \mapsto g \circ f$ is a $\quant$-functor in each variable separately, since $\hom(g,g') \preceq \hom(gf, g'f)$ and $\hom(f,f') \preceq \hom(gf, gf')$.
It is not, in general, a $\quant$-functor out of the product $[\D,\E] \times [\C,\D]$ of \Cref{ex:constructions}(b), which would require $\hom(g,g') \meet \hom(f,f') \preceq \hom(gf, g'f')$.
Over $\law$, let $\C = \D = \E$ be $\{0, 1, \dots, 10\}$ with $\hom(x,y) = \max(0, y-x)$, let $f = g = \id$, and let $f'(x) = g'(x) = \min(x+1, 10)$. Then $\hom(g,g') \meet \hom(f,f')$ is the cost $1$, whereas $\hom(gf, g'f')$ is the cost $2$.
When $\quant$ is commutative, composition is a $\quant$-functor out of the tensor product, whose hom-objects are the products $\hom(g,g') \odot \hom(f,f')$, and $\QCat$ is enriched in the closed symmetric monoidal category of $\quant$-categories~\citep[Sections~1.4 and~2.3]{Kelly2005}.

The following proposition is the universal property of the weighted global sections.
For a $\quant$-category $\C$ and a cellular $\quant$-sheaf $\shf$, write $[\C, \shf(-)] : \cat{G} \to \QCat$ for the functor $\sigma \mapsto [\C, \shf(\sigma)]$, which acts on incidences by postcomposition.

\begin{prop}\label{prop:weighted-limit}
    Let $\shf$ be a cellular $\quant$-sheaf and $W$ a weight function on $G$. For every $\quant$-category $\C$, there is an isomorphism of $\quant$-categories
    \[ \bigl[\C, \wsections\bigr] \cong \bigl[\widetilde{W}, [\C,\shf(-)]\bigr], \]
    natural in $\C$, which sends a $\quant$-functor $\Lambda$ to the natural transformation $\beta$ with $\beta_\sigma(a)(c) = \Lambda(c)_\sigma(a)$ for every object $a$ of $\widetilde{W}(\sigma)$.
\end{prop}

\begin{proof}
    For a family of maps $\lambda_\sigma : \obj(\C) \times \obj(\widetilde{W}(\sigma)) \to \obj(\shf(\sigma))$ indexed by $\sigma \in V \sqcup E$, set $\Lambda(c)_\sigma(a) = \beta_\sigma(a)(c) = \lambda_\sigma(c,a)$. At a vertex $v \in V$, since $\widetilde{W}(v) = \terminal$, the component $\beta_v : \terminal \to [\C,\shf(v)]$ is a $\quant$-functor if and only if its single value $x_v \coloneqq \beta_v(\ast) : \C \to \shf(v)$ is a $\quant$-functor. At an edge $e \in E$ with $\partial e = \{u,v\}$, each $\beta_e(u), \beta_e(v) : \C \to \shf(e)$ is a $\quant$-functor if and only if $\hom_\C(c,c') \preceq \hom_{\shf(e)}(\lambda_e(c,u), \lambda_e(c',u))$ and similarly at $v$, and $\beta_e : \widetilde{W}(e) \to [\C,\shf(e)]$ is a $\quant$-functor if and only if $W(u,v) \preceq \hom_{[\C,\shf(e)]}(\beta_e(u), \beta_e(v)) = \bigmeet_{c \in \C} \hom_{\shf(e)}(\lambda_e(c,u), \lambda_e(c,v))$ and $W(v,u) \preceq \bigmeet_{c \in \C} \hom_{\shf(e)}(\lambda_e(c,v), \lambda_e(c,u))$. Naturality of $\beta$ along the incidences $u \face e$ and $v \face e$ holds if and only if $\beta_e(u) = \shf_{u \face e} \circ x_u$ and $\beta_e(v) = \shf_{v \face e} \circ x_v$. By \Cref{thm:global-sections}, these conditions hold if and only if $\Lambda(c) = (x_v(c))_{v \in V}$ belongs to $\wsections$ for every $c \in \C$ and $\Lambda : \C \to \wsections$ is a $\quant$-functor. For two such families $\lambda$ and $\lambda'$, with associated $\quant$-functors $\Lambda$ and $\Lambda'$ and natural transformations $\beta$ and $\beta'$, the hom-objects $\hom(\Lambda,\Lambda')$ and $\hom(\beta,\beta')$ are both equal to $\bigmeet_{c,\sigma,a} \hom_{\shf(\sigma)}(\lambda_\sigma(c,a), \lambda'_\sigma(c,a))$, where $a$ ranges over the objects of $\widetilde{W}(\sigma)$, and naturality in $\C$ follows by precomposition.
\end{proof}

Recall from \Cref{sec:global-sections} that, over every unital quantale $\quant$, the $\quant$-category $\wsections$ is the end $\int_{\sigma \in \cat{G}} [\widetilde{W}(\sigma), \shf(\sigma)]$.
When $\quant$ is commutative, \Cref{prop:weighted-limit} is the defining universal property of the limit of $\shf$ weighted by $\widetilde{W}$ in the $\QCat$-enriched sense~\citep[Section~3.1]{Kelly2005}, with $\QCat$ enriched in $(\QCat,\tensor,\terminal)$ and $\cat{G}$ regarded as a free $\QCat$-category.
Under the correspondence between ordinary functors on $\cat{G}$ and $\QCat$-functors on the free $\QCat$-category on $\cat{G}$~\citep[Section~2.5]{Kelly2005}, the end above is the hom-object of the functor category, and naturality in $\C$ is $\QCat$-naturality because two $\quant$-functors that agree on objects are equal.
Over a noncommutative base we take \Cref{prop:weighted-limit} as the universal property of the weighted global sections, which it determines up to isomorphism, as the case $\C = \terminal$ shows.

\bibliographystyle{plainnat}
\bibliography{biblio}

\end{document}